\documentclass[12pt,a4paper]{amsproc}

\usepackage{amsmath}
\usepackage{amsfonts}
\usepackage{amssymb}

\usepackage[dvips,colorlinks,bookmarks,breaklinks]{hyperref} 
\usepackage{graphicx}       
\usepackage[T1]{fontenc}
\usepackage{lmodern}        

\newcommand{\change}[1][]{}
\newcommand{\bem}[1]{}

\allowdisplaybreaks[1]

\newtheorem{theorem}{Theorem}[section]
\newtheorem{lemma}[theorem]{Lemma}
\newtheorem{proposition}[theorem]{Proposition}

\theoremstyle{definition}
\newtheorem{definition}[theorem]{Definition}

\theoremstyle{remark}
\newtheorem{remark}[theorem]{Remark}

\numberwithin{equation}{subsection}
\numberwithin{theorem}{subsection}

\newcommand{\A}{\ensuremath{\mathbb{A}}}
\newcommand{\C}{\ensuremath{\mathbb{C}}}

\newcommand{\R}{\ensuremath{\mathbb{R}}}
\newcommand{\Z}{\ensuremath{\mathbb{Z}}}
\newcommand{\N}{\ensuremath{\mathbb{N}}}
\newcommand{\PP}{\ensuremath{\mathbb{P}^1_\R}}
\newcommand{\eps}{\ensuremath{\varepsilon}}
\newcommand{\HH}{\ensuremath{\mathbb{H}}}
\newcommand{\id}{\ensuremath{\mathbf{1}}}

\newcommand{\cA}{\ensuremath{\mathcal{A}}}
\newcommand{\cB}{\ensuremath{\mathcal{B}}}

\newcommand{\cJ}{\ensuremath{\mathcal{J}}}
\newcommand{\cL}{\ensuremath{\mathcal{L}}}
\newcommand{\cM}{\ensuremath{\mathcal{M}}}
\newcommand{\cN}{\ensuremath{\mathcal{N}}}

\newcommand{\rR}{\ensuremath{\mathrm{R}}}

\newcommand{\TT}[1][]{\ensuremath{{T^\prime_{#1}}}}
\newcommand{\dd}{\ensuremath{\mathrm{d}}}

\newcommand{\lb}{\ensuremath{[\![}}             
\newcommand{\rb}{\ensuremath{]\!]}}             
\newcommand{\lbd}{\ensuremath{[\![}}            
\newcommand{\rbd}{\ensuremath{]\!]^\star}}      
\newcommand{\lbs}{\ensuremath{\big(}}           
\newcommand{\rbs}{\ensuremath{\big)}}           

\newcommand{\Iq}{\ensuremath{I_q}}
\newcommand{\Ir}{\ensuremath{I_{R_q}}}
\newcommand{\Iqs}{\ensuremath{I_q^\star}}
\newcommand{\Irs}{\ensuremath{I_{R_q}^\star}}

\newcommand{\simGn}{\ensuremath{\sim_\Gn{q}}}
\newcommand{\simR}{\ensuremath{\sim_{\text{reg}}}}

\newcommand{\Ar}{\ensuremath{\cA_q^\text{reg}}}   
\newcommand{\Adr}{\ensuremath{\cA_q^\text{dreg}}} 

\newcommand{\re}[1]{\ensuremath{{\mathrm{Re}\left(#1\right)}}}
\newcommand{\im}[1]{\ensuremath{{\mathrm{Im}\left(#1\right)}}}
\newcommand{\sign}[1]{\ensuremath{{\mathrm{sign}\left(#1\right)}}}
\newcommand{\Gn}[1]{\ensuremath{{G_{#1}}}}
\newcommand{\SL}[1]{\ensuremath{{\mathrm{SL}\!\left(2, #1 \right)}}}
\newcommand{\PSL}[1]{\ensuremath{{\mathrm{PSL}\!\left(2, #1 \right)}}}
\newcommand{\rMatrix}[4]{{\begin{bmatrix} #1 & #2 \\ #3 & #4 \end{bmatrix}}}
\newcommand{\Matrix}[4]{{\begin{pmatrix} #1 & #2 \\ #3 & #4 \end{pmatrix}}}

\newcommand{\abs}[1]{\ensuremath{{\left\lvert#1\right\rvert}}}
\newcommand{\ov}[1]{\ensuremath{\overline{ #1 }}}
\newcommand{\nextinteger}[1]{\ensuremath{\left\lfloor #1\right\rfloor}}
\newcommand{\nli}[1]{\ensuremath{\left\langle #1 \right\rangle_q}}
\newcommand{\nlid}[1]{\ensuremath{\left\langle #1 \right\rangle_q^\star}}
\newcommand{\orbit}[1]{\ensuremath{\mathrm{orbit}\left( #1 \right)}}  
\newcommand{\sorbit}[1]{\ensuremath{\mathrm{orbit}^\star\left( #1 \right)}}  

\title{Nearest $\lambda_q$-multiple fractions}

\author{D.~Mayer}
\address{\href{http://www.dynamik.tu-clausthal.de}{Institute of Theoretical Physics}, \href{http://www.tu-clausthal.de}{TU Clausthal}, 38678- Clausthal-Zellerfeld, Germany}
\email{\href{mailto:dieter.mayer@tu-clausthal.de}{dieter.mayer@tu-clausthal.de}}

\author{T.~M\"uhlenbruch}
\address{\href{http://www.fernuni-hagen.de/WTHEORIE/tobias.muehlenbruch.html}{Department of  Mathematics and Computer Science}, \href{http://www.fernuni-hagen.de}{Fern\-Universit\"at in Hagen}, 58084- Hagen, Germany}
\email{\href{mailto:Tobias.Muehlenbruch@fernuni-hagen.de}{Tobias.Muehlenbruch@fernuni-hagen.de}}
\thanks{The second named author was supported in part by the Deutsche Forschungsgemeinschaft through the DFG Research Project ``Transfer operators and non arithmetic quantum chaos'' (Ma 633/16-1).}


\date{\today}


\begin{document}


\begin{abstract}
We discuss the nearest $\lambda_q$--multiple continued fractions and their duals for $\lambda_q=2 \cos\left( \frac{\pi}{q}\right)$ which are closely related to the Hecke triangle groups $\Gn{q}$, $q=3,4,\ldots$.
They have been introduced in the case $q=3$ by Hurwitz and for even $q$ by Nakada.
These continued fractions are generated by interval maps $f_q$ respectively $f_q^\star$ which are conjugate to subshifts over infinite alphabets.
We generalize to arbitrary $q$ a result of Hurwitz concerning the $\Gn{q}$- and $f_q$-equivalence of points on the real line.
The natural extension of the maps $f_q$ and $f_q^\star$ can be used as a Poincar\'e map for the geodesic flow on the Hecke surfaces $\Gn{q} \backslash \HH$ and allows to construct the transfer operator for this flow.
\end{abstract}


\maketitle



\setcounter{tocdepth}{3} \tableofcontents



\section{Introduction}
\label{A}
In the transfer operator approach to Selberg's zeta-function for Fuchsian groups $G$ \cite{Ma91}, \cite{Ma} this function is expressed through the Fredholm-determinant of the generalized Perron-Frobenius operator $\cL_\beta$ for the geodesic flow on the corresponding surface $G\backslash \HH$ of constant negative curvature.
This operator is constructed through an expanding interval map $f \colon I \to I$ closely related to a Poincar\'e map of the flow.
In all the cases treated up to now this interval map generates some kind of continued fraction expansion like the Gauss expansion or its extensions such that the length spectrum of the flow can be completely characterized by the periodic orbits of $f$ respectively the purely periodic continued fraction expansions.
This program has been carried out in full detail for the modular surfaces $\Gamma\setminus \HH$ defined by subgroups $\Gamma \subset \PSL{\Z}$ of the full modular group.

For these groups the transfer operator has another rather important property:
its eigenfunctions with eigenvalue $1$ can be directly related to their automorphic forms, that is real analytic Eisenstein series and Maass wave forms respectively the holomorphic modular forms.
This relation gave rise to the theory of periodic functions \cite{LZ01}, \cite{BLZ} which generalize the Eichler-Manin-Shimura cohomology theory for holomorphic modular forms.

The physical interpretation of these relations between the transfer operator and the spectral properties of the Laplacian for these groups $G$ is within the theory of quantum chaos \cite{Ru08}, \cite{Sa95}:
the transfer operator encodes the classical length spectrum of the geodesic flow and relates these data to the quantum data, namely eigenvalues and eigenfunctions respectively resonances of its quantized system.
In this sense this transfer operator approach extends the more common approach to quantum chaos via the Selberg-Gutzwiller trace formula \cite[Theorem~13.8, p.~209]{Hej83}, \cite{Marklof03}.

Obviously it is necessary to work out the transfer operator for more general Fuchsian groups, especially non-arithmetic ones,
for which the Hecke triangle groups $\Gn{q}$ are good examples, since up to the cases $q=3$, $4$, $6$ all of them are indeed non-arithmetic.
In \cite{MS08} the authors constructed a symbolic dynamics for the geodesic flow on the Hecke surfaces $\Gn{q}\backslash \HH$ for arbitrary $q$, the case $q=3$ was treated earlier in \cite{KU07}.
In both cases the authors used the nearest $\lambda_q$-multiple continued fraction expansion, denoted for short by $\lambda_q$-CF, and its dual expansion.
Another approach was discussed also in \cite{SS95}.
Some of the ergodic properties of these  $\lambda_q$-CF's for $q$ even have been worked out in \cite{Na95} by H.\ Nakada.
In the present paper we discuss the $\lambda_q$-CF's and their duals for arbitrary $q$ via their generating interval maps $f_q$ and $f_q^\star$, which allow us to derive also a transfer operator for the Hecke triangle groups $\Gn{q}$, whose Fredholm determinant is closely related to the Selberg function for the groups $\Gn{q}$ as we will discuss in a forthcoming paper.

\medskip

In \cite{Hu89} Hurwitz introduced nearest integer continued fraction expansions of the form
\numberwithin{equation}{section}
\begin{equation}
\label{A.1}
a_0 + \frac{-1 }{ a_1 + \frac{-1}{a_2 + \frac{-1}{a_3 + \frac{-1}{\ldots}}}}
\end{equation}
where $a_0$ is an arbitrary integer and the $a_i$, $i \geq 1$, are integers satisfying $\abs{a_i} \geq 2$ and $a_i \, a_{i+1} <0$ if $\abs{a_i}=2$.
They are generated by the interval map
\numberwithin{equation}{subsection} 
\[
f_3 \colon I_3  \to I_3; \quad  x \mapsto \frac{-1}{x} - \left\langle \frac{-1}{x} \right\rangle,
\]
where $I_3 = \left[ -\frac{1}{2}, \frac{1}{2}\right]$ and $\langle x \rangle$ denotes the nearest integer to $x$,
by the usual algorithm:
\begin{itemize}
\item[(0)] $a_0 = \left\langle x\right\rangle$ and $x_1:= x-a_0$,
\item[(1)] $a_1 = \left\langle \frac{-1}{x_1}\right\rangle$ and $x_2:= \frac{-1}{x_1} - a_1 = f_3(x_1)$,
\item[(\,$i$\,)] $a_i = \left\langle \frac{-1}{x_i}\right\rangle$ and $x_{i+1}:= \frac{-1}{x_i} - a_i = f_3(x_i)$.
\item[($\star$)] The algorithm terminates if $x_{i+1}=0$.
\end{itemize}

Let $\PSL{\Z} = \SL{\Z} \bmod \{\pm \id\}$ denote the full modular group.
Elements of the group can be identified with $2 \times 2$-matrices with integer entries and determinant $1$, up to a common sign.
The group acts on the projective real line $\R\cup\{\infty\}$ by M\"obius transformations ${a \; b \choose c \; d} \, z = \frac{az + b}{cz + d}$.
The group $\PSL{\Z}$ is generated by the elements $S$ and $T$ corresponding to the actions $z \mapsto \frac{-1}{z}$ and $z \mapsto z+1$.
The generators satisfy the relations $S^2 = (ST)^3 = \id$.
In particular, the elements $T^a$ and $ST^a$ correspond to the actions $z \mapsto z+a$ and $z \mapsto \frac{-1}{a+z}$.
Hence we can write the continued fraction expansion in \eqref{A.1} in terms of a (formal) M\"obius transformation as $T^{a_0} \, ST^{a_1} \, ST^{a_2} \, ST^{a_3} \, \cdots \; 0$.

Hurwitz found in \cite{Hu89}, that equivalence of two points $x,y \in \R$ under the generating map $f_3$ and its extension to $\R$ is not the same as equivalence under the group action of $\PSL{\Z}$.
This is obviously in contrast with the case of the Gauss map $f_\text{G} \colon[0,1]\to [0,1]$ with $f_G(x)=\frac{1}{x}\mod 1$ and the modular group $\PSL{\Z}$.

\smallskip

In \cite{Na95} Nakada introduced for even integers $q \geq 4$ the nearest $\lambda_q$-multiple continued fractions with $ \lambda_q = 2 \cos \frac{\pi}{q}$, which we will denote by $\lambda_q$-CF's.
They are similar to the Rosen contined fractions introduced in \cite{Ro54} and discussed in detail in \cite{BKS99}.
The extension to the case $q\geq 3$ odd is straightforward, where $q=3$ corresponds to the nearest integer continued fractions of Hurwitz.
These $\lambda_q$-CF's and their dual expansions, introduced for $q=3$ also by Hurwitz, can be generated by interval maps $f_q$ and $f_q^\star$ closely related to the Hecke triangle groups $\Gn{q}$.
Both maps are conjugate to subshifts over infinite alphabets, which when reduced to certain sofic systems, determine completely the properties of the corresponding $\lambda_q$-CF and its dual expansion.
It turns out, that Hurwitz's result on equivalence of points on the real axis under $f_3$ and the group action of $\PSL{\Z}$ is true for general $q \geq 3$:
there exists for every $q\geq 3$ exactly one pair of points $\big(r_q,-r_q\big)$ which are equivalent under $\Gn{q}$ but not under the map $f_q$.
The natural extension $F_q$ of the interval map $f_q \colon \Iq \to \Iq$ can be easily constructed from the symbolic dynamics of the maps $f_q$ and $f_q^\star$ as sofic systems.
It can be used to construct a Poincare section for the geodesic flow on the Hecke surface $\Gn{q} \backslash \HH$, and hence also a transfer operator for the group $\Gn{q}$ and its Selberg zeta function.
The properties of this operator will be discussed elsewere.

\medskip

The structure of this article is as follows:
In Section~\ref{B} we introduce the Hecke triangle groups and the $\lambda_q$-CF's respectively the dual $\lambda_q$-CF's.
In Section~\ref{C} we discuss the interval maps $f_q$ and $f_q^\star$ generating the nearest $\lambda_q$-multiple continued fractions and construct Markov partitions for these maps.
In Section~\ref{D} we show that the maps $f_q$ and $f_q^\star$ are conjugate to subshifts over infinite alphabets and introduce sofic systems closely related to the $\lambda_q$-CF and its dual.
This allows a simple construction of the natural extension $F_q$ of the map $f_q$.
In Section~\ref{E} we relate the natural extension $F_q$ to the geodesic flow on the Hecke surfaces
$\Gn{q}\backslash \HH$ and derive the transfer operator for this flow.
The final Section~\ref{F} contains a discussion of the convergence properties of the $\lambda_q$-CF's by relating them to reduced Rosen $\lambda$-fractions as discussed in \cite{Ro54}.

\section{Nearest $\lambda_q$-multiple continued fractions}
\label{B}

\subsection{Hecke triangle groups}
\label{B1}
Hecke triangle groups are Fuchsian \linebreak groups of the first kind, all except three are non-arithmetic.
The reader may wish to consult~\cite[pp.\ 591, 627]{He83} for a discussion of Hecke triangle groups and their relation to Dirichlet series.

Denote by $\PSL{\R}$ the projective special linear group given by
\begin{equation}
\label{B1.1}
\PSL{\R} = \SL{\R} / \left\{ \pm \id \right\}
\end{equation}
where $\pm \id = \Matrix{\pm 1}{0}{0}{\pm 1}$.
We denote by $\rMatrix{a}{b}{c}{d} = \left\{ \Matrix{a}{b}{c}{d}, \Matrix{-a}{-b}{-c}{-d}\right\}$ the elements of $\PSL{\R}$, but identify often elements of $\PSL{\R}$ and $\SL{\R}$.

For a given integer $q \geq 3$ the $q^\mathrm{th}$ \emph{Hecke triangle group} $\Gn{q}$ is generated by
\begin{equation}
\label{B1.2}
S:= \rMatrix{0}{-1}{1}{0} \quad \text{and} \quad T_q := \rMatrix{1}{\lambda_q}{0}{1}
\end{equation}
with relations
\begin{equation}
\label{B1.3}
S^2 = \left( ST_q \right)^q = \id,
\end{equation}
where $\lambda_q$ is given by
\begin{equation}
\label{B1.4}
\lambda_q := 2\cos\left( \frac{\pi}{q} \right).
\end{equation}
Later on we also need the element
\begin{equation}
\label{B1.5}
\TT[q] := \rMatrix{1}{0}{\lambda_q}{1} = ST_q^{-1}S \in \Gn{q}.
\end{equation}
We may suppress the $q$-dependence in our notation when we work with a fixed value of $q$.

The Hecke triangle group $\Gn{q}$ is a discrete subgroup of $\PSL{\R}$ and its limit set is the projective line $\PP = \R \cup \{\infty\}$.
It acts on the upper half-plane, the lower half-plane and on $\PP$ by \emph{M\"obius transformations}
\begin{equation}
\label{B1.6}
\Gn{q} \times \PP \to \PP; \quad \left(\rMatrix{a}{b}{c}{d} , x  \right) \mapsto \rMatrix{a}{b}{c}{d} \, x :=
\begin{cases}
\frac{ax + b}{cx + d}  & \text{ if } x \in \R \mbox{ and} \\
\frac{a}{c}            & \text{ if } x = \infty.
\end{cases}
\end{equation}

\medskip

The points $x,y \in \PP$ are \emph{$\Gn{q}$-equivalent} denoted by $x \simGn y$, if there exists an element $g \in \Gn{q}$ such that $y= g\,x$.
Obviously, this is an equivalence relation.

\subsection{Nearest $\lambda_q$-multiple continued fractions and their duals}
\label{B2}
Consider finite or infinite sequences $\lbs a_i\rbs_i$.
We denote periodic parts of the sequences by overlining the period part and finitely repeated patterns are denoted by a power where a $0^\text{th}$ power vanishes:
\begin{align*}
\lbs a_1, \ov{a_2, a_3} \rbs
  &= \lbs a_1, \,a_2, a_3,\,a_2, a_3,\,a_2, a_3,\,\ldots \rbs, \\
\lbs a_1, (a_2, a_3)^i,a_4, \ldots \rbs
&= \lbs a_1, \underbrace{a_2, a_3,a_2, a_3, \ldots, a_2,a_3}_{i \, \mathrm{times} \; a_2,a_3}, a_4, \ldots \rbs \quad \text{and} \\
\lbs a_1, (a_2)^0, a_3, \ldots \rbs
  &= \lbs a_1, a_3, \ldots \rbs.
\end{align*}
We use also $- \lbs a_1,\ldots \rbs = \lbs -a_1,\ldots \rbs$.

Put
\begin{equation}
\label{B2.1}
h_q :=
\begin{cases}
\frac{q-2}{2} \quad & \text{for $q$ even and} \\
\frac{q-3}{2} \quad & \text{for $q$ odd.}
\end{cases}
\end{equation}
We define the set $\cB_q$ of \emph{forbidden blocks} as
\begin{equation}
\label{B2.3}
\cB_q:=
\begin{cases}
\{ \lbs \pm 1 \rbs\} \cup \bigcup_{m=1}^\infty \{ \lbs\pm 2, \pm m\rbs\}
  & \text{for $q=3$}, \\
\{ \lbs (\pm 1)^{h_q+1} \rbs\} \cup  \bigcup_{m=1}^\infty \{ \lbs(\pm 1)^h_q, \pm m\rbs\}
  & \text{for $q$ even and} \\
\{ \lbs (\pm 1)^{h_q+1} \rbs\} \cup \\
\; \cup \bigcup_{m=1}^\infty \{ \lbs(\pm 1)^{h_q},\pm 2,(\pm 1)^{h_q},\pm m\rbs\}
  & \text{for $q$ odd, $q \geq 5$.}
\end{cases}
\end{equation}
The choice of the sign is the same within each block.
For example $\lbs 2,3 \rbs$, $\lbs -2,-3 \rbs \in \cB$ and $\lbs 2,-3 \rbs \not\in \cB$ for $q=3$.

We call a sequence $\lbs a_1,a_2,a_3,\ldots \rbs$ \emph{$q$-regular} if $\lbs a_l,a_{l+1}, \ldots, a_L \rbs \not\in \cB_q$ for all $1 \leq l < L$ and \emph{dual $q$-regular} if $\lbs a_L,a_{L-1}, \ldots, a_l \rbs \not\in \cB_q$ for all $1 \leq l < L$.
Denote by $\Ar$ respectively by $\Adr$ the set of infinite $q$-regular respectively dual $q$-regular sequences $(a_i)_{i\in\N}$.
\smallskip

A \emph{nearest $\lambda_q$-multiple continued fraction}, or $\lambda_q$-CF, is a formal expansion of the type
\begin{eqnarray}
\label{B2.4}
[a_0; a_1,a_2,a_3,\ldots]:= a_0 \lambda_q + \frac{-1}{a_1\lambda_q + \frac{-1}{a_2\lambda_q + \frac{-1}{a_3\lambda_q + \ldots } } }
\end{eqnarray}
with $a_i \in \Z_{\not=0}$, $i \geq 1$ and $a_0 \in \Z$.

We say that $[a_0; a_1,a_2,a_3,\ldots]$ \emph{converges} if either $[a_0; a_1,a_2,a_3,\ldots] = [a_0; a_1,a_2,a_3,\ldots,a_L]$ has finite length or $\lim_{L\to \infty} [a_0; a_1,a_2,a_3,\ldots, a_L]$ exists in $\R$.

We adopt the same notations as  introduced for sequences earlier.
For example we write $[ a_0; a_1, \ov{a_2, a_3}]$ for a periodic tail of the expansion and $- [ a_0; a_1,\ldots ]$ for $[ -a_0; -a_1,\ldots ]$.

A $\lambda_q$-CF is \emph{regular} respectively \emph{dual regular} if the sequence \linebreak $\lbs a_1,a_2,a_3,\ldots \rbs$ is $q$-regular respectively dual $q$-regular.
Regular and dual regular $\lambda_q$-CF's are denoted by $\lb a_0; a_1,\ldots \rb$ respectively $\lbd a_0; a_1,\ldots \rbd$.

\begin{proposition}
\label{B2.5}
Regular and dual regular $\lambda_q$-CF's converge.
\end{proposition}

\begin{proof}
This follows immediately from Lemmas~4 and~34 in \cite{MS08}.
\end{proof}

An alternative proof of Proposition~\ref{B2.5} for infinite regular and dual regular expansions with leading $0$ follows also from the lemmas in the Sections~\ref{D2a} and \ref{D3a}.

\smallskip

Converging $\lambda_q$-CF's can be rewritten in terms of elements of the Hecke triangle group $\Gn{q}$:
if the expansion~\eqref{B2.4} is finite it can be written as follows
\begin{align}
\nonumber
[a_0; a_1, a_2, a_3, \ldots, a_L]
&= a_0 \lambda_q + \frac{-1}{a_1\lambda_q + \frac{-1}{a_2\lambda_q + \frac{-1}{a_3\lambda_q + \ldots \frac{-1}{a_L\lambda_q} } }} \\
\label{B2.6}
&= T^{a_0} \, ST^{a_1} \, ST^{a_2} \, ST^{a_3} \, \cdots \, ST^{a_L}\; 0,
\end{align}
since $\frac{-1}{a\lambda_q + x} = ST^a \; x$.
For infinite converging $\lambda_q$-CF the expansion has to be interpreted as
\begin{align*}
[a_0; a_1, a_2, a_3, \ldots]
&= \lim_{L \to \infty}[a_0; a_1, a_2, a_3, \ldots, a_L] \\
&= \lim_{L\to \infty} T^{a_0} \, ST^{a_1} \, ST^{a_2} \, ST^{a_3} \, \cdots \, ST^{a_L}\; 0 \\
&= T^{a_0} \, ST^{a_1} \, ST^{a_2} \, ST^{a_3} \, \cdots \; 0.
\end{align*}

An immediate consequence of this is
\begin{lemma}
\label{B2.7}
For a finite regular $\lambda_q$-CF one finds for $q$ even
\begin{align*}
\lb a_0; a_1, \ldots,a_n,(1)^h \rb    &= \lb a_0; a_1, \ldots,a_n-1,(-1)^h \rb \\
 \quad \text{respectively for $q$ odd} \\
\lb a_0; \ldots,a_n,(1)^h,2,(1)^h \rb &= \lb a_0; \ldots,a_n-1,(-1)^h,-2,(-1)^h \rb.
\end{align*}
\end{lemma}

\begin{proof}
Assume the left hand side to be regular.
This implies $a_n \not= 1$ and hence the right hand side is regular, too.
Conversely, assume the right hand side to be regular and hence $a_n-1 \not= -1$.
Therefore the expansions on the left hand side are regular.

The identity now follows by writing $\lambda_q$-CF's in terms of M\"obius transformations and using the identity $(ST)^{h_q} \; 0 = T^{-1}\,(ST^{-1})^{h_q} \; S T^{-1}S \; 0  = T^{-1}\,(ST^{-1})^{h_q} \; 0$ since $0$  is a fixed point of $S T^{-1} S  = \TT$.
\end{proof}

\smallskip

\begin{remark}
\label{B2.8}
For $q=3$ the nearest $\lambda_q$-multiple continued fractions are in fact the well-known nearest integer fractions extensively discussed by Hurwitz in \cite{Hu89}.
In particular, Theorem~\ref{B5.2} for $q=3$ was proved by him there.
We include his results for the sake of completeness and to show how this special case $q=3$ fits well into the discussion of the case of odd $q \geq 5$.
See also Remark~\ref{B3.9}.
\end{remark}

\begin{remark}
\label{B2.9}
For $q \geq 4$ the regular $\lambda_q$-CF's correspond to Rosen's \linebreak $\lambda_q$-fractions introduced in \cite{Ro54} and discussed in \cite{BKS99}.
We will discuss this relation in more detail in \S\ref{F1}.
\end{remark}

\subsection{Special values and their expansions}
\label{B3}
The following results are shown in \cite{MS08}:

The point $x=\mp\frac{\lambda_q}{2}$ has the regular $\lambda_q$-CF
\begin{equation}
\label{B3.2}
\mp\frac{\lambda_q}{2} = \begin{cases}
\lb 0;(\pm 1)^{h_q} \rb         \qquad & \mbox{for even } q \mbox{ and} \\
\lb 0;(\pm 1)^{h_q},\pm 2,(\pm 1)^{h_q} \rb \qquad & \mbox{for odd } q.
\end{cases}
\end{equation}

Put
\begin{equation}
\label{B3.8}
\begin{split}
R_q &:= \lambda_q + r_q
\quad \text{with} \\
r_q &:=
\begin{cases}
\lb 0; \ov{3} \rb & \text{for $q = 3$,} \\
\lb 0; \ov{(1)^{h_q-1},2} \rb & \text{for $q$ even and} \\
\lb 0; \ov{(1)^{h_q},2,(1)^{h_q-1},2} \rb & \text{for $q$ odd, $q \geq 5$}.
\end{cases}
\end{split}
\end{equation}
whose expansion is periodic of length $\kappa_q$ with
\begin{equation}
\label{B3.1}
\kappa_q :=
\begin{cases}
h_q = \frac{q-2}{2}& \text{for even $q$ and} \\
2h_q+1 = q-2& \mbox{for odd } q,
\end{cases}
\end{equation}
The regular respectively dual regular $\lambda_q$-CF of the point $x=R_q$ has the form
\begin{align}
\label{B3.3}
\qquad
R_q &=
\begin{cases}
~ \lb 1 ; \ov{(1)^{h_q-1},2} \rb          \quad & \mbox{for even } q, \\
~ \lb 1 ; \ov{(1)^{h_q},2,(1)^{h_q-1},2} \rb \quad & \mbox{for odd } q \geq 5 \mbox{ and}\\
~ \lb 1 ; \ov{3} \rb                     \quad & \mbox{for } q=3.
\end{cases} \\
\label{B3.4}
&=
\begin{cases}
\lbd 0; (-1)^h,\ov{-2,(-1)^{h_q-1}} \rbd                 & \text{for even } q, \\
\lbd 0; (-1)^h,\ov{-2,(-1)^{h_q},-2,(-1)^{h_q-1}} \rbd   & \text{for odd } q \geq 5 \text{ and}\\
\lbd 0; -2,\ov{-3} \rbd                                  & \text{for } q=3.
\end{cases}
\end{align}
Moreover,
\begin{align}
\label{B3.5}
&R_q = 1  \quad \text{and} \quad  -R_q = S \, R_q
&& \text{for even } q \text{ and} \\
\label{B3.5a}
&R_q^2 +(2-\lambda_q)R_q = 1 \quad \text{and} \quad  -R_q = \big(TS\big)^{h_q+1} \, R_q
&& \text{for odd }  q
\end{align}
and $R_q$ satisfies the inequality
\begin{equation}
\label{B3.6}
\frac{\lambda_q}{2} < R_q \leq 1.
\end{equation}

\begin{remark}
\label{B3.7}
For $R_3$ one finds
\[
1+R_3 = \frac{1+\sqrt{5}}{2}.
\]
\end{remark}

\begin{remark}
\label{B3.9}
The form of the $\lambda_q$-CF of $r_3$ in~\eqref{B3.8} can be obtained from the expansions for $q$ odd, $q\geq 5$ by interpreting it as a M\"obius transformation with $(1)^{-1}$ as $ST^{-1}$:
\begin{align*}
r_3
&=
\lb 0; \ov{1^{h_3},2,(1)^{h_3-1},2} \rb = \lb 0; \ov{2,(1)^{-1},2} \rb \\
&= ST^2 \, ST^{-1} \, ST^2 \,\cdot\, ST^2 \, ST^{-1} \, ST^2 \; \cdots \; 0 \\
&= ST^2 \, TSTS \, ST^2 \,\cdot\, ST^2 \, TSTS \, ST^2 \; \cdots \; 0 \\
&= ST^3 \, ST^3 \,\cdot\, ST^3 \, ST^3 \; \cdots \; 0 = \lb 0; \ov{3} \rb.
\end{align*}
\end{remark}

\subsection{A lexiographic order}
\label{B4}
Let $x,y \in \Ir :=\left[ -R_q,R_q \right]$ have the regular $\lambda_q$-CF's $x = \lb a_0; a_1, \ldots \rb$ and $y = \lb b_0; b_1, \ldots \rb$.
Denote by $l(x)\leq \infty$ respectively $l(y)\leq \infty$ the number of entries in the above $\lambda_q$-CF's.
We introduce a \emph{lexiographic order} ``$\prec$'' for $\lambda_q$-CF's in the following way:
For $n \in \Z_{\geq 0}$ being the number of identical digits at the head of the $\lambda_q$-CF's, i.e., $a_i =b_i$ for all $0 \leq i \leq n$ and $l(x), l(y) \geq n$, we define
\begin{equation}
\label{B4.2}
x \prec y :\iff
\begin{cases}
a_0 < b_0     & \text{if $n=0$}, \\
a_n > 0 > b_n & \text{if $n>0$, both $l(x), l(y) \geq n+1$ and $a_n b_n < 0$}, \\
a_n < b_n     & \text{if $n>0$, both $l(x), l(y) \geq n+1$ and $a_n b_n > 0$}, \\
b_n < 0       & \text{if $n>0$ and $l(x)=n$ or} \\
a_n > 0       & \text{if $n>0$ and $l(y)=n$}.
\end{cases}
\end{equation}
We also write $x \preceq y$ for $x \prec y$ or $x = y$.

This is indeed an order on regular $\lambda_q$-CF's, since Lemmas~22 and~23 in \cite{MS08} imply:
\begin{lemma}
\label{B4.3}
Let $x$ and $y$ have regular $\lambda_q$-CF's.
Then $x \prec y$ if and only if $x < y$.
\end{lemma}

\smallskip

The authors of \cite{MS08} introduce a process called ``rewriting'' of $\lambda_q$-CF's where   forbidden blocks in the $\lambda_q$-CF are replaced by allowed ones without changing its value.
The rules for ``rewriting'' are based on the interpretation of a $\lambda_q$-CF in terms of M\"obius transformations given by  group elements of the Hecke group, see~\eqref{B2.6}, and applying the group relations~\eqref{B1.3}.
We refer  in particular to Lemma~11 and Lemma~13 in \cite{MS08} for the details.
A simple example for this rewriting is used in the proof of Lemma~\ref{B2.7}.

It follows  from the proof of Lemma~34 in \cite{MS08} that every dual regular $\lambda_q$-CF can be rewritten into a regular $\lambda_q$-CF.

\begin{lemma}
\label{B4.4}
The lexiographic order $\prec$ in \eqref{B4.2} can be extended to dual regular $\lambda_q$-CF's with leading digit $0$.
Rewriting two dual regular $\lambda_q$-CF's satisfying $ \lbd 0;a_1,\ldots \rbd \prec \lbd 0; b_1,\ldots \rbd$ into regular $\lambda_q$-CF's does not change their order.
\end{lemma}

\begin{remark}
\label{B4.1}
The lexicographic order $"\prec"$ however cannot be defined for all dual regular $\lambda_q$-CF's with arbitrary leading coefficient as the following example shows:
consider the dual regular $\lambda_3$-CF's of $R_3$ in \eqref{B3.3} and \eqref{B3.4}.
Obviously $R_3 = \lbd 0; -2,\ov{-3} \rbd  =   \lbd 1 ; \ov{3} \rbd$.
Extending naively ``$\prec$'' in \eqref{B4.2} to this case would lead to $\lbd 0; -2,\ov{-3} \rbd  \prec  \lbd 1 ; \ov{3} \rbd$ and hence $\lbd 0; -2,\ov{-3} \rbd \prec \lbd 0; -2,\ov{-3} \rbd$.
\end{remark}

\begin{proof}[Proof of Lemma~\ref{B4.4}]
The $\lambda_q$-CF's $\lbd 0;a_1,\ldots \rbd$ and $\lbd 0; b_1,\ldots \rbd$ are dual regular.
No rewriting is necessary if both are also regular.

Assume $\lbd 0;b_1,b_2,\ldots \rbd$ starts with a forbidden block.
If it is of the form
\[
\lbd 0;b_1,b_2,\ldots \rbd =
\begin{cases}
\lbd 0;(1)^{h_q},m \rbd               \quad& \text{for even $q$,} \\
\lbd 0;(1)^{h_q},2,(1)^{h_q},m \rbd   \quad& \text{for odd $q\geq 5$ and} \\
\lbd 0;2,m \rbd                       \quad& \text{for $q=3$}
\end{cases}
\]
with $m \geq 2$ for $q\geq 4$ respectively $m \geq 3$ for $q=3$, then $\lbd 0;a_1,a_2,\ldots \rbd$ must be of the form
\[
\lbd 0;a_1,a_2,\ldots \rbd =
\begin{cases}
\lbd 0;(1)^{h_q},n \rbd               \quad& \text{for even $q$,} \\
\lbd 0;(1)^{h_q},2,(1)^{h_q},n \rbd   \quad& \text{for odd $q\geq 5$ and} \\
\lbd 0;2,n \rbd                       \quad& \text{for $q=3$}
\end{cases}
\]
with $n<m$ and $n \geq 2$ for $q\geq 4$ respectively $n \geq 3$ for $q=3$.
Using the rewriting rules in Lemmas~11 and~13 in \cite{MS08} we find
\begin{align*}
\lbd 0;a_1,a_2,\ldots \rbd \!\!
&\to
\lb \tilde{a}_0;\tilde{a}_1,\tilde{a}_2,\ldots \rb \\
&:=
 \begin{cases}
\lb -1;(-1)^{h_q},n-1,\ldots \rb                \;& \text{for even $q$,} \\
\lb -1;(-1)^{h_q},-2,(-1)^{h_q},n-1,\ldots \rb  \;& \text{for odd $q\geq 5$,} \\
\lb -1;-2,n-1,\ldots \rb                        \;& \text{for $q=3$ and}
\end{cases}
\end{align*}
respectively
\begin{align*}
\lbd 0;b_1,b_2,\ldots \rbd \!\!
&\to
\lb \tilde{b}_0;\tilde{b}_1,\tilde{b}_2,\ldots \rb \\
&:=
 \begin{cases}
\lb -1;(-1)^{h_q},m-1,\ldots \rb                \;& \text{for even $q$,} \\
\lb -1;(-1)^{h_q},-2,(-1)^{h_q},m-1,\ldots \rb  \;& \text{for odd $q\geq 5$,} \\
\lb -1;-2,m-1,\ldots \rb                        \;& \text{for $q=3$}
\end{cases}
\end{align*}
hence by \eqref{B4.2} $\lb -1;\tilde{a}_1,\tilde{a}_2,\ldots\rb \prec \lb -1; \tilde{b}_1,\tilde{b}_2,\ldots \rb$.

If $\lbd 0;b_1,b_2,b_3,\ldots\rbd$ is of the form
\begin{equation}
\label{B4.5}
\lbd 0;b_1,b_2,\ldots \rbd
= \begin{cases}
~\lbd 0;(-1)^{h_q},-m , \ldots \rbd           \quad& \text{for even $q$,} \\
~\lbd 0;(-1)^{h_q},-2,(-1)^{h_q},-m, \ldots \rbd   & \text{for odd $q\geq 5$,} \\
~\lbd 0;(-2,-m), \ldots \rbd                  \quad& \text{for $q=3$}
\end{cases}
\end{equation}
with $m \geq 2$ for $q \geq 4$ respectively $m \geq 3$ for $q=3$, and $\lbd 0;a_1,a_2,\ldots \rbd$ does not contain a forbidden block starting with $a_1$, the rewriting rules in \cite{MS08} give
\begin{eqnarray*}
\lbd 0;b_1,b_2,\ldots \rbd
&\to&
\lb \tilde{b}_0;\tilde{b}_1,\tilde{b}_2,\ldots \rb \\
&=&
 \begin{cases}
\lb 1;(1)^{h_q},1-m,\ldots \rb                \;& \text{for even $q$,} \\
\lb 1;(1)^{h_q},2,(1)^{h_q},1-m,\ldots \rb    \;& \text{for odd $q\geq 5$,} \\
\lb 1;2,1-m,\ldots \rb                        \;& \text{for $q=3$.}
\end{cases}
\end{eqnarray*}
Therefore \eqref{B4.2} implies $\lbd 0;a_1,a_2,\ldots\rbd  \prec \lb 1; \tilde{b}_1,\tilde{b}_2,\ldots \rb$.

If $\lbd 0;b_1,b_2,b_3,\ldots \rbd$ is of the form~\eqref{B4.5} and $\lbd 0;a_1,a_2,\ldots \rbd$ is of the form
\[
\lbd 0;a_1,a_2,\ldots \rbd =
\begin{cases}
\lbd 0;(-1)^{h_q},-n \rbd              \quad& \text{for even $q$,} \\
\lbd 0;(-1)^{h_q},2,(1)^{h_q},-n \rbd  \quad& \text{for odd $q\geq 5$ and} \\
\lbd 0;2,-n \rbd                       \quad& \text{for $q=3$}
\end{cases}
\]
with $n>m$ then the rewriting rules in \cite{MS08} lead to
\begin{eqnarray*}
&& \lb 1;\tilde{a}_1,\tilde{a}_2,\ldots \rb =
 \begin{cases}
\lb 1;(1)^{h_q},1-n,\ldots \rb                \;& \text{for even $q$,} \\
\lb 1;(1)^{h_q},2,(1)^{h_q},1-n,\ldots \rb    \;& \text{for odd $q\geq 5$,} \\
\lb 1;2,1-n,\ldots\rb                         \;& \text{for $q=3$,}
\end{cases} \\
&&
\lb 1;\tilde{b}_1,\tilde{b}_2,\ldots \rb =
 \begin{cases}
\lb 1;(1)^{h_q},1-m,\ldots \rb              \;& \text{for even $q$,} \\
\lb 1;(1)^{h_q},2,(1)^{h_q},1-m,\ldots \rb  \;& \text{for odd $q\geq 5$,} \\
\lb 1;2,1-m,\ldots \rb                      \;& \text{for $q=3$,}
\end{cases}
\end{eqnarray*}
and hence $\lb 1;\tilde{a}_1,\tilde{a}_2,\ldots \rb \prec\lb 1;\tilde{b}_1,\tilde{b}_2,\ldots \rb$.

Completely analogous are the cases when $\lbd 0;a_1,\ldots \rbd$ starts with a forbidden block or the first forbidden block starts at $a_n,\,n>1$ and $\lbd 0;b_1,b_2,\ldots \rbd$ is a regular $\lambda_q$-CF.
If both $\lbd 0;a_1,\ldots \rbd$ and $\lbd 0;b_1,\ldots \rbd$ have the same forbidden block starting at $a_1$ respectively $b_1$ then both dual regular $\lambda_q$-CF's are rewritten in the same way and the forbidden block does not influence the order ``$\prec$''.
\end{proof}

\subsection{Equivalence relations and continued fractions}
\label{B5}

Let $x,y \in \R$ have infinite regular $\lambda_q$-CF's $x = \lb a_0; a_1, \ldots \rb$ and $y = \lb b_0; b_1, \ldots \rb$.
We say that $x$ and $y$ are \emph{regular $\lambda_q$-CF-equivalent}, denoted by $x \simR y$, if the regular $\lambda_q$-CF's of $x$ and $y$ have the same tail, i.e., there exists $m,n \in \N$ such that the sequences $(a_m, a_{m+1}, \ldots)$ and $(b_n, b_{n+1}, \ldots)$ coincide.
Obviously, this is an equivalence relation.
We can extend this equivalence relation to all regular $\lambda_q$-CF's by declaring all finite regular $\lambda_q$-CF's to be regular $\lambda_q$-CF-equivalent.

\begin{theorem}[Equivalence relations]
\label{B5.2}
For $x,y\in \mathbb{R}$ the following properties are equivalent:
\begin{enumerate}
\item[(1)]
$x \simGn y$.
\item[(2)]
$x$ and $y$ satisfy:
\begin{itemize}
\item $x \simR y$ or
\item $x \simR \pm r$ and $y \simR \mp r$.
\end{itemize}
\end{enumerate}
\end{theorem}

To prove the proposition, we need the following lemmas:

\begin{lemma}
\label{B5.6}
If $x$ has an infinite regular $\lambda_q$-CF and $g \in \Gn{q}$ satisfies $g\,x \in \R$, then $g \, x$ has an infinite $\lambda_q$-CF with at most $h_q$ consecutive digits $\pm 1$.
Its tail coincides with the tail of the regular $\lambda_q$-CF of $x$.
\end{lemma}

\begin{proof}
Let $x$ have the regular $\lambda_q$-CF $x= \lb a_0; a_1, \ldots \rb$.
We can write $g$ as a word in the generators $S$ and $T_q$: $g = T_q^{b_0} \, ST_q^{b_1} \, ST_q^{b_2} \, \cdots \, ST_q^{b_m}\, S^\delta$ with $b_0 \in \Z$, $b_i \in \Z_{\not=0}$, $i=1,\ldots,m$, and $\delta \in \{0,1\}$.
Then $g\,x$ can formally be written as $g \,x = T^{b_0} \, ST^{b_1} \, ST^{b_2} \, \cdots \, ST^{b_m}\, S^\delta \; T^{a_0} \, ST^{a_1} \, ST^{a_2} \, \cdots \; 0$.

Consider for $n > m$ sufficiently large the element $g_n\in \Gn{q}$ given by
\begin{equation}
\label{B5.7}
g_n := T_q^{b_0} \, ST_q^{b_1} \, ST_q^{b_2} \, \cdots \, ST_q^{b_m}\, S^\delta \; T_q^{a_0} \, ST_q^{a_1} \, ST_q^{a_2} \, \cdots \, ST_q^{a_n} \in \Gn{q}.
\end{equation}
The identities $S^2 = \id$, $\left( ST_q^{\pm1} \right)^q = \id$, $T_q^a \, \left(ST_q^{\pm1} \right)^{q-1} \,ST_q^b = T_q^{a+b \mp 1}$ and \linebreak $T_q^a \, \left(ST_q^{\pm1} ST_q^b \right)^l = T_q^{a \mp 1}\, \left(ST_q^{\mp 1} \right)^{q-l-2} \,ST_q^{b \mp 1}$ for $h_q+1 \leq l \leq q-2$ and $a,b$ arbitrary follow from~\eqref{B1.3}.
But $q-l-2 \leq q-(h_q+1)-2 = h_q-1$ for $q$ even and $q-l-2 \leq h_q$ for $q$ odd.
We apply these identities recursively on $g_n$ in \eqref{B5.7}.
After a finite number of steps one arrives at a word representing $g_n $ which contains blocks of at most $h_q$ consecutive digits $\pm 1$.
Indeed, since each application of one of these identities reduces the length of the word, the process of applying the identities has to stop after a finite number of steps.
And, since the $\lambda_q$-CF of $x$ is reduced, there are no blocks of more than $h_q$ consecutive $\pm 1$ to the right of the right of part ``$ST^{b_m}\, S^\delta \; T^{a_0} \, ST^{a_1}$'' of the word $g_n$ in \eqref{B5.7}.

Hence $g\, x$ can be written as a $\lambda_q$-CF of the form~\eqref{B2.6} without blocks of more than $h_q$ consecutive digits $\pm 1$ and with a tail identical to the regular tail in the $\lambda_q$-CF of $x$.
\end{proof}

\begin{lemma}
\label{B5.4}
Let $\big[a_0;a_1,a_2,\ldots\big]$ be an infinite $\lambda_q$-CF containing no blocks of more than $h_q+1$ consecutive digits $\pm 1$ and at most one block of $h_q+1$ consecutive digits $\pm 1$ for $q \geq 4$ respectively no forbidden digits $\pm 1$ for $q =3$.
If the block $\big[(\pm1)^{h_q+1}\big]$ exists the block has to be the first forbidden block of the $\lambda_q$-CF and has to be preceded by a digit of alternate sign.

If the first forbidden block starts at $a_i$, $i \geq 1$ and its rewriting leads to a new forbidden block then this forbidden block and its rewritten version must have the form given in Table~\ref{B5.5}.
The new forbidden block will appear to the right of $a_i$.
If the new forbidden block is of the form $\big[ (\pm1)^{h_q+1} \big]$ then its preceding digit is negative.
\end{lemma}

\begin{table}
\begin{center}
\begin{footnotesize}
\begin{tabular}{|l c l| }
\hline
even $q$: ($m \geq 1$)&& \\
$\pm \big[a_{i-1},\underline{1^{h_q+1}},2,1^{h_q-1},m\big]$
& $\!\rightarrow\!$ &
$\pm \big[a_{i-1}-1,(-1)^{h_q-1},\underline{(1)^{h_q},m}\big]$ \\
$\pm \big[a_{i-1},\underline{1^{h_q},2},1^{h_q-1},m\big]$
& $\!\rightarrow\!$ &
$\pm \big[a_{i-1}-1,(-1)^{h_q},\underline{1^{h_q},m}\big]$ \\
\hline
odd $q \geq 5$: ($m \geq 2$)&& \\
$\pm \big[a_{i-1},\underline{1^{h_q+1}},2,1^{h_q} \big]$
& $\!\rightarrow\!$ &
$\pm \big[a_{i-1}-1,(-1)^{h_q},\underline{1^{h_q+1}}\big]$ \\
$\pm \big[a_{i-1},\underline{1^{h_q+1}},2,1^{h_q-1},2,1^{h_q},m\big]$
& $\!\rightarrow\!$ &
$\pm \big[a_{i-1}-1,(-1)^{h_q},\underline{1^{h_q},2,1^{h_q},m}\big]$ \\
$\pm \big[a_{i-1},\underline{1^{h_q},2,1^{h_q},2},1^{h_q}\big]$
& $\!\rightarrow\!$ &
$\pm \big[a_{i-1}-1,(-1)^{h_q},-2,(-1)^{h_q},\underline{1^{h_q+1}}\big]$ \\
$\pm \big[a_{i-1},\underline{1^{h_q},2,1^{h_q},2},1^{h_q-1},2,1^{h_q},m \big]$
& $\!\rightarrow\!$ &
$\pm \big[a_{i-1}-1,(-1)^{h_q},-2,(-1)^{h_q},\underline{1^{h_q},2,1^{h_q},m}\big]$ \\
\hline
$q=3$: ($m \geq 3$ and $n \geq 0$)&&\\
$\pm \big[a_{i-1},\underline{2,3},m\big]$
& $\!\rightarrow\!$ &
$\pm \big[ a_{i-1}-1,-2,\underline{2,m}\big]$ \\
$\pm \big[a_{i-1},\underline{2,2},2^n,3,m\big]$
& $\!\rightarrow\!$ &
$\pm \big[ a_{i-1}-1,-(3+n),\underline{2,m}\big]$ \\
\hline
\end{tabular}
\end{footnotesize}
\end{center}

\caption{
Under assumptions of Lemma~\ref{B5.4} we list all possibilities where the rewriting of a forbidden block generates a new forbidden block.
The forbidden blocks are underlined.
}
\label{B5.5}
\end{table}

\begin{proof}
W.l.o.g.\ assume that the forbidden block starting at $a_i$ has positive digits and hence $a_{i-1} \not= 1$.
For $q$ even the forbidden block must have the form $\big[ 1^{h_q+1} \big]$ or $\big[1^{h_q},m \big]$ with $m \geq 2$.
The rewriting rules in Lemma~11 of \cite{MS08} lead to
\begin{align*}
\big[\ldots,a_{i-1}, 1^{h_q+1},a_{i+h_q+1},\ldots \big]
&\to \big[\ldots,a_{i-1}-1, (-1)^{h_q-1}, a_{i+h_q+1}-1,\ldots \big] \text{ and}\\
\big[\ldots,a_{i-1}, 1^{h_q}, m,a_{i+h_q+1},\ldots \big]
&\to \big[\ldots,a_{i-1}-1, (-1)^{h_q}, m-1 , a_{i+h_q+1},\ldots \big]
\end{align*}
for $m \geq 2$.
Changing $a_{i-1}$ to $a_{i-1}-1$ cannot introduce a forbidden block since $a_{i-1} \not=1$ and a digit $a_{i-1} \geq 2$ cannot follow a block of the from $\big[ 1^{h_q} \big]$.
Hence any new forbidden block has to start with the digit $a_{i+h_q+1}-1$ respectively $m-1$.
Two cases are possible:
$\big[ a_{i+h_q+1}, \ldots \big] = \big[2,1^{h_q-1},l,\ldots]$ respectively $\big[m,\ldots \big]=\big[2,1^{h_q-1},l,\ldots]$ with $l \geq 1$.
This shows that the block $\big[a_i, \ldots \big]$ must have the form $\big[ 1^{h_q+1},2,1^{h_q-1},l\big]$ respectively $\big[ 1^{h_q},2,1^{h_q-1},l \big]$ with its rewriting leading to the form as stated in the lemma.

For $q$ odd, $q \geq 5$, the forbidden block $[a_i,\ldots]$ has either the form $\big[1^{h_q},2,1^{h_q},m\big]$ with $m \geq 2$ or the form $\big[1^{h_q+1}\big]$.
Rewriting rules in Lemma~13 of \cite{MS08} then give
\begin{align*}
\big[\ldots,a_{i-1}, 1^{h_q+1},a_{i+h_q+1},\ldots \big]
&\to \big[\ldots,a_{i-1}-1, (-1)^{h_q}, a_{i+h_q+1}-1,\ldots \big] \text{ and}\\
\big[\ldots,a_{i-1}, 1^{h_q}, 2,1^{h_q}, m,\ldots \big]
&\to \big[\ldots,a_{i-1}-1, (-1)^{h_q},-2,(-1)^{h_q}, m-1,\ldots \big]
\end{align*}
for $m \geq 2$, and similar arguments as for $q$ even show that the forbidden block and the digits following it are either of the form $\big[ 1^{h_q+1} \big]$ followed by $\big[2,1^{h_q}\big]$ or $\big[ 2,1^{h_q-1},2,1^{h_q},l \big]$, $l \geq 2$, or $\big[ 1^{h_q},2,1^{h_2},2 \big]$ followed by $\big[ 1^{h_q} \big]$ respectively $\big[ 1^{h_q-1},2,1^{h_q},l \big]$, $l \geq 2$.
The rewritten form is the as given in the lemma.
For $a_{i-1}=2$ rewriting cannot lead to a new forbidden block to the left of $a_i$ contradicting otherwise the first forbidden block to start with $a_i$.

The case $q=3$ with forbidden block $\big[ 2,2^n,m \big],\,m \geq 2$ and $n \in \Z_{\geq 0}$ can be handled in complete analogy by using the rewriting rule $\big[ a,2,2^n,b \big] \to \big[a-1,-2-n,b-1 \big]$ with $a,b \not=2$.
\end{proof}

\begin{proof}[Proof of Proposition~\ref{B5.2}]
We show first the implication $(2) \Rightarrow (1)$.
If $x \simR y$ then $x$ and $y$ have regular $\lambda_q$-CF's with the same tail:
\[
x=\lb a_0;a_1,\ldots,a_m,a_{m+1},\ldots \rb
\text{ and }
y=\lb b_0;b_1,\ldots,b_n,a_{m+1},\ldots \rb.
\]
Put $g:=  T^{a_0} \, ST^{a_1} \, \cdots \, ST^{a_m}   \left( T^{b_0} \, ST^{b_1} \, \cdots \, ST^{b_n} \right)^{-1} \in \Gn{q}$.
Writing $x$ and $y$ in terms of  M\"obius transformations as explained in \eqref{B2.6} we find,
\begin{eqnarray*}
g \, y
&=&
g \, T^{b_0}\, ST^{b_1} \, \cdots \, ST^{b_n} \, ST^{a_{m+1}} \, \cdots \; 0 \\
&=&
T^{a_0}\, ST^{a_1} \, \cdots \, ST^{a_m} \, ST^{a_{m+1}} \, \cdots \; 0 = x
\end{eqnarray*}
and hence $x$ and $y$ are $\Gn{q}$-equivalent.

Assume next $x \simR r$ and $y \simR -r$ and hence $x \simGn r$ and $y \simGn -r$.
Since $R_q=T_q\, r_q $ according to  \eqref{B3.8} and $-R_q=S R_q$ according to \eqref{B3.5} for even $q$ respectively $-R_q=(T_q)^{h_q+1} R_q $ according to \eqref{B3.5a} for odd $q$ obviously $r \simGn -r$ and hence $x \simGn y$.

\smallskip

To show implication $(1) \Rightarrow (2)$, assume there exists $g \in \Gn{q}$ with $g\, x = y$ with $x$ and $y$ having infinite regular $\lambda_q$-CF's.
Lemma~\ref{B5.6} shows that $g \, x$ can be written as an infinite $\lambda_q$-CF with regular tail satisfying the assumptions of Lemma~\ref{B5.4}.
Using the rewriting rules in Lemmas~11 and~13 of \cite{MS08} we can recursively rewrite the $\lambda_q$-CF of $g\,x$ into a regular $\lambda_q$-CF from the left to the right.
We procede to the next forbidden block if rewriting does not lead to a new forbidden block.
Lemma~\ref{B5.4} implies that a new forbidden block can only appear to the right of the original one which we process next.
If this rewriting process stops after finitely many steps then $y=g\,x$ and $x$ have the same tail in their $\lambda_q$-CF's and $x \simR y$.

Hence assume, the rewriting process has to be repeated again and again.
Then after a sufficiently large but finite number of rewriting steps one arrives at the situation where the $\lambda_q$-CF of $g\,x$ is regular up to one forbidden block.
Denote this $\lambda_q$-CF by $[a_0;a_1,\ldots]$ with the remaining forbidden block starting at digit $a_i$, $i \geq 1$ and assume w.l.o.g.\ the forbidden block has positive digits.

Consider first the case $q$ even:
By Lemma~\ref{B5.4} the forbidden block and the following digits have the form $[a_i,\ldots] = [B_0,B_1,B_2,\ldots]$ with the block $B_0 \in \left\{ \big[ 1^{h_q+1},2 \big], \big[ 1^{h_q},2 \big]\right\}$ and the blocks $B_j \in \left\{ \big[ 1^{h_q},2 \big], \big[ 1^{h_q-1},2 \big]\right\}$, for all $j \geq 1$.
Since by assumption $B_0$ was the last forbidden block in the $\lambda_q$-CF, necessarily $B_j\neq \big[ 1^{h_q},2 \big]$.
Hence the $\lambda_q$-CF of $g\, x$ has the form
\[
g\, x = \big[a_0;a_1,\ldots, a_{i-1}, \underline{1^{l},2}, \ov{1^{h_q-1},2}  \big]
\quad \text{with } l = h_q, h_q+1
\]
where the forbidden block at digit $a_i$ is underlined, and whose tail, determining also the tail of $x$, is  regular $\lambda_q$-CF-equivalent to $r_q$.
After infinitely many further rewritings one arrives at the regular $\lambda_q$-CF of $y$ whose tail is regular $\lambda_q$-CF-equivalent to $-r_q$.

Consider next the case $q \geq 5$ odd:
Lemma~\ref{B5.4} again determines the form of the forbidden block and the following digits as
\[
[a_i,\ldots] = [B_0,B_1,B_2,\ldots]
\]
with the block $B_0 \in \left\{ \big[ 1^{h_q+1},2 \big], \big[ 1^{h_q},2,1^{h_q},2 \big]\right\}$ and the blocks
\[
B_j \in \left\{ A_1:=\big[ 1^{h_q},2 \big], A_2:=\big[ 1^{h_q-1},2,1^{h_q},2 \big]\right\},
\qquad j \geq 1.
\]
Since the blocks $\big[A_1, A_1\big]$ and $\big[A_2, A_1\big]$ are forbidden blocks, necessarily $B_j = A_2$ for all $j \geq 2$, since
otherwise $B_0$ would not be the last forbidden block in the $\lambda_q$-CF of $g\,x$.
Hence the $\lambda_q$-CF of $g\, x$ has the form
\[
g\, x = \big[a_0;a_1,\ldots, a_{i-1}, \underline{B_0}, B_1, \ov{1^{h_q-1},2,1^{h_q},2}  \big]
\]
where the forbidden block at digit $a_i$ is again underlined.
As in the previous case, we find $x$ is regular $\lambda_q$-CF-equivalent to $r_q$ and $y$ is regular $\lambda_q$-CF-equivalent to $-r_q$.

Consider finally the case $q = 3$:
Lemma~\ref{B5.4} gives again the form of the forbidden block and the following digits as $[a_i,\ldots] = [B_0,B_1,\ldots]$ with  $B_0 = \big[ 2,2^n,3 \big]$, $n \geq 0$, and the blocks $B_j \in \left\{ \big[ 2 \big], \big[ 3 \big]\right\}$, $j \geq 1$.
Since the blocks $[2,2]$ and $[2,3]$ are forbidden, necessarily $B_j = [3]$ for all $j \geq 1$.
The $\lambda_q$-CF of $g\, x$ hence has the form
\[
g\, x = \big[a_0;a_1,\ldots, a_{i-1}, \underline{2,2^n,3}, \ov{3}  \big]
\]
where the forbidden block at digit $a_i$ is underlined.
Again $x$ is regular $\lambda_q$-CF-equivalent to $r_q$ and $y$ is regular $\lambda_q$-CF-equivalent to $-r_q$.

\end{proof}

\section{Generating maps for the $\lambda_q$-continued fractions and their duals}
\label{C}

Similar to the Gauss continued fractions also the $\lambda_q$-continued fractions and their duals, which
for $q=3$ have been introduced by Hurwitz in \cite{Hu89}, can be generated by interval maps with strong
ergodic properties like in the case of the Gauss maps.

\subsection{The interval maps $f_q$ and $f_q^\star$}
\label{C1}
Denote by $\Iq$ respectively $\Ir$ the intervals
\begin{equation}
\label{C1.1}
\Iq = \left[ -\frac{\lambda_q}{2}, \frac{\lambda_q}{2} \right]
\quad \text{respectively} \quad
\Ir = \left[-R_q,R_q\right]
\end{equation}
with $\lambda_q$ as in \eqref{B1.4} and $R_q= \lambda_q + r_q$ as in \eqref{B3.8}.
The \emph{nearest $\lambda_q$-multiple map} $\nli{\cdot}$ is given by
\begin{equation}
\label{C1.2}
\nli{\cdot} \colon \R \to \Z; \quad x \mapsto \nli{x}:= \nextinteger{ \frac{x}{\lambda_q} +
\frac{1}{2} }
\end{equation}
where $\nextinteger{\cdot}$ is the (modified) floor function
\begin{equation}
\label{C1.3}
\nextinteger{x} = n
\iff
\begin{cases}
n < x \leq n+1 \quad& \text{if } x > 0 \text{ and} \\
n \leq x < n+1 \quad& \text{if } x \leq 0.
\end{cases}
\end{equation}

We also need the map $\nlid{\cdot}$ given by
\begin{equation}
\label{C1.4}
\nlid{\cdot} \colon \R \to \Z; \quad x \mapsto \nlid{x} :=
\begin{cases}
\nextinteger{ \frac{x}{\lambda_q} + 1 - \frac{R_q}{\lambda_q} }   \qquad & \text{if } x  \geq 0 \text{ and} \\
\nextinteger{ \frac{x}{\lambda_q} + \frac{R_q}{\lambda_q} }       \qquad & \text{if } x < 0.
\end{cases}
\end{equation}

\smallskip

The interval maps $f_q \colon \Iq \to \Iq$ and $f_q^\star \colon \Ir \to \Ir$ are defined as follows:
\begin{align}
\label{C1.5}
f_q(x) &= \begin{cases}
-\frac{1}{x} - \nli{\frac{-1}{x}} \lambda_q  \qquad & \text{if } x \in \Iq \backslash \{0\},\\
0 & \text{if } x =0 \,\,\text{and}
\end{cases} \\
\label{C1.6}
f_q^\star(y) &= \begin{cases}
-\frac{1}{y} - \nlid{\frac{-1}{y}} \lambda_q  \qquad & \text{if } y \in \Ir \backslash \{0\},\\
0 & \text{if } y =0.
\end{cases}
\end{align}

\subsection{$\lambda_q$-CF's and their generating interval maps}
\label{C2}
The maps $f_q$ and $f_q^\star$ generate the regular respectively dual regular $\lambda_q$-CF's in the following sense:

For given $x, y \in \R$ the coefficients $a_i$ and $b_1$, $i \in \Z_{\geq 0}$ are determined by the following algorithms:
\begin{itemize}
\item[(0)] $a_0 = \nli{x}$ and $x_1:= x-a_0\lambda_q \in \Iq$,
\item[(1)] $a_1 = \nli{\frac{-1}{x_1}}$ and $x_2:= \frac{-1}{x_1}-a_1\lambda_q = f_q(x_1) \in \Iq$,
\item[(\,$i$\,)] $a_i = \nli{\frac{-1}{x_i}}$ and $x_{i+1}:= \frac{-1}{x_i}-a_i\lambda_q = f_q(x_i) \in \Iq$,
\item[($\star$)] The algorithm terminates if $x_{i+1}=0$
\end{itemize}
and
\begin{itemize}
\item[(0)] $b_0 = \nlid{y}$ and $y_1:= y-b_0\lambda_q\in \Ir$,
\item[(1)] $b_1 = \nlid{\frac{-1}{y_1}}$ and $y_2:= \frac{-1}{y_1}-b_1\lambda_q= f_q^\star(y_1) \in \Ir$,
\item[(\,$i$\,)] $b_i = \nlid{\frac{-1}{y_i}}$ and $y_{i+1}:= \frac{-1}{y_i}-b_i\lambda_q= f_q^\star(y_i) \in \Ir$,
\item[($\star$)] The algorithm terminates if $y_{i+1}=0$.
\end{itemize}
By construction the coefficients form $\lambda_q$-CF's in the sense of~\eqref{B2.4}:
\begin{equation}
\label{C2.1}
x = [a_0; a_1, a_2, \ldots ]
\quad \text{and} \quad
y = [b_0; b_1, b_2, \ldots ].
\end{equation}

\begin{proposition}
\label{C2.2}
The $\lambda_q$-CF of $x$ in \eqref{C2.1} is unique for all $x$ not in $\bigcup_{n=0}^\infty f_q^{-n}(\pm \frac{\lambda_q}{2})$ and regular whereas the one of $y$ is unique for all $y\notin \bigcup_{n=0}^\infty \left(f_q^\star\right)^{-n} (\pm r_q)$ and dual regular.
\end{proposition}

\begin{proof}
A simple calculation shows that the regular $\lambda_q$-CF of all points $x=\pm \frac{2}{2m-1)\lambda_q},\; m=2,3,\ldots$ and their preimages is not unique. But these points belong to the preimages of the points $\pm \frac{\lambda_q}{2}$  On the other hand the dual  $\lambda_q$-CF of the points $y=\pm\frac{1}{r_q+m\lambda_q},\; m=1,2,\ldots$ and their preimages is not unique. But these points are all the preimages of the points $\pm r_q$.
\end{proof}

\begin{remark}
The non-uniqueness of certain finite regular $\lambda_q$-CF's in Lemma~\ref{B2.7} can also be derived from Proposition~\ref{C2.2}.
\end{remark}

\subsection{Markov partitions for $f_q$ and $f_q^\star$}
\label{C3}
Obviously $f_q$ is locally expanding, that means $\abs{f_q^\prime(x)} > 1$ for all $x \in \Iq$, if one takes the one-sided derivatives at the points of discontinuity.
The same holds true for the map $f_q^\star$ for $q$ odd.
For $q$ even ${f_q^\star}^\prime( \pm R_q) = 1$ but $\abs{\left(f_q^{\star2}\right)^\prime(y)} > 1$ for all $y \in \Ir$, and hence both maps $f_q$ and $f_q^\star$ are locally smooth, expanding maps.
Indeed both maps have the Markov property, that means that they allow for Markov partitions.
To construct these partitions we use the orbits of the boundary points of the two intervals $\Iq$ and $\Ir$ respectively the monotonicity intervals of the maps $f_q$ and $f_q^\star$.

Define the orbit of $x$ under $f_q$ respectively $f_q^\star$ as
\begin{align}
\label{C3.1}
\orbit{x}  &= \left\{ x, f_q(x), f_q^2(x) :=  f_q(f_q(x)), f_q^3(x), \ldots  \right\}  \\
\nonumber
           &= \left\{ f_q^n(x);\; n =0,1,2,\ldots \right\}  \quad \text{respectively} \\
\sorbit{x} &= \left\{ \left(f_q^\star\right)^n(x);\; n =0,1,2,\ldots \right\}.
\end{align}
The orbits $\orbit{-\frac{\lambda_q}{2}}$ and $\sorbit{-R_q}$ are both finite.
If $\sharp \{S\}$ denotes the cardinality of the set $S$, we have
\[
\sharp \orbit{-\frac{\lambda_q}{2}} = \sharp \sorbit{-R_q} = {\kappa_q} +1,
\]
as can be seen from the regular $\lambda_q$-CF of $-\frac{\lambda_q}{2}$ in~\eqref{B3.2} and the dual regular $\lambda_q$-CF of $-R_q$ in~\eqref{B3.4}.
We denote the elements of $\orbit{-\frac{\lambda_q}{2}}$ by $\phi_i$ respectively of $\sorbit{-R_q}$ by $\psi_i$, $i=0,\ldots, {\kappa_q}$, such that
\begin{align}
\label{C3.2}
-R_q = -\psi_{0} &< -\frac{\lambda_q}{2} = \phi_0 < \psi_{1} < \phi_1 < \psi_{2} < \phi_2 < \ldots\\
\nonumber
& \qquad
\ldots < \psi_{{\kappa_q}-2} < \phi_{{\kappa_q}-2} < \psi_ {{\kappa_q}-1}< \phi_{{\kappa_q}-1} = - \frac{1}{\lambda_q} < \psi_{{\kappa_q}} < \phi_{{\kappa_q}} = 0
\end{align}
holds.
By using the regular $\lambda_q$-CF of $-\frac{\lambda_q}{2}$ and the dual regular $-\lambda_q$-CF of $R_q$
respectively the order ``$\prec$'' in \S\ref{B4} one easily verifies
\begin{lemma}
\label{C3.3}
The order in \eqref{C3.2} is achieved for $q$ even by defining
\begin{equation}
\label{C3.4}
\phi_i = f_q^i \left( - \frac{\lambda_q}{2} \right)
\quad\text{and} \quad
\psi_i = \left(f_q^\star\right)^i (-R_q),\quad
0 \leq i \leq h_q= {\kappa_q},
\end{equation}
respectively for $q$ odd by defining
\begin{align}
\label{C3.5}
\phi_{2i} &= f_q^i \left( - \frac{\lambda_q}{2} \right),
\qquad
\phi_{2i+1} = f_q^{h+i+1} \left( - \frac{\lambda_q}{2} \right)
\quad\text{and}\\
\nonumber
\psi_{2i} &= \left(f_q^\star\right)^i (-R_q),
\qquad
\psi_{2i+1} = \left(f_q^\star\right)^{h+i+1} (-R_q),
\quad 0 \leq i \leq h_q=\frac{{\kappa_q}-1}{2}.
\end{align}
\end{lemma}

In the case $q=3$ one has $\kappa_3=1$ and $h_3=0$.
Therefore
\[
\phi_0 = -\frac{1}{2}, \quad \phi_1 = 0, \quad \psi_0 = -R_3 = \frac{1-\sqrt{5}}{2} \text{ and } \psi_1 = R_3-1 = \frac{\sqrt{5}-3}{2}.
\]

\smallskip

Define next $\phi_{-i}=-\phi_i$, $0 \leq i \leq {\kappa_q}$, respectively $\psi_{-i}= -\psi_i$ for $0 \leq i \leq {\kappa_q+1}$ with $\psi_{\kappa_q+1}=0$.

Obviously the intervals
\begin{equation}
\label{C3.6}
\Phi_i   := \big[ \phi_{i-1},\phi_i \big]
\quad \text{and} \quad
\Phi_{-i}:= \big[ \phi_{-i} ,\phi_{-(i-1)} \big]
 \qquad 1 \leq i \leq {\kappa_q}
\end{equation}
respectively
\begin{equation}
\label{C3.7}
\Psi_i    := \left[ \psi_{i-1},\psi_i \right]
\quad \text{and} \quad
\Psi_{-i} := \left[ \psi_{-i} ,\psi_{-(i-1)} \right],
\end{equation}
$1 \leq i \leq {\kappa_q+1}$, define Markov partitions of the intervals $\Iq$ and $\Ir$:
this means that
\begin{align*}
&\bigcup_{\eps =+,-} \bigcup_{i=1}^{\kappa_q} \Phi_{\eps i} = \Iq,
\quad
 \Phi_{\eps i}^\circ \cap \Phi_{\delta j}^\circ= \emptyset
\quad \text{for} \quad \eps i \neq \delta j \\
&\bigcup_{\eps = +,-} \bigcup_{i=1}^{{\kappa_q+1}} \Psi_{\eps i} = \Ir,
\quad
\Psi_{\eps i}^\circ \cap \Psi_{\delta j}^\circ = \emptyset
\quad \text{for} \quad \eps i \neq \delta j .
\end{align*}
where $S^\circ$ denotes the interior of the set $S$.
To get a reasonable symbolic dynamics for the two maps $f_q$ and $f_q^\star$ we have to construct finer
partitions using the monotonicity intervals of the two maps.
\medskip
Consider first the case $q=3$ such that $\lambda_3=1$.
Define for $m=2,3,4,\ldots$ the intervals $J_m$ as
\begin{equation}
\label{C3.8}
J_2 = \left[ -\frac{1}{2}, -\frac{2}{5} \right]
\quad \text{and} \quad
J_m = \left[ -\frac{2}{2m-1}, -\frac{2}{2m+1} \right], \quad m = 3,4,\ldots,
\end{equation}
and set $J_{-m}:= -J_m$ for $m=2,3,4,\ldots$.
Since $f_3\left(J_{\pm 2}\right) = \mp \left[0,\frac{1}{2}\right]$ and $f_3(J_{\pm m}) = I_3$ for $m=3,4,\ldots$
the partition satisfies
\[
\bigcup_{\eps = +,-}\bigcup_{m=2}^\infty J_{\eps m} = I_3
\quad \text{and} \quad
 J_{\eps m}^\circ \cap J_{\delta k}^\circ = \emptyset \quad \text{for} \quad\eps m \neq \delta k.
\]
Hence this partition, which we denote by $\cM (f_3)$, is Markovian.
The maps $f_3 \big|_{J_m}$ are monotone with $f_3 \big|_{J_m}(x) = -\frac{1}{x}-m$ and locally invertible with $\left(f_3 \big|_{J_m}\right)^{-1}(y) = -\frac{1}{y+m}$ for $y \in f_3(J_m)$.

For $q \geq 4$ define intervals $J_m$, $m=1,2,\ldots$, as
\begin{equation}
\label{C3.9}
\begin{split}
J_1 &= \left[ -\frac{\lambda_q}{2}, -\frac{2}{3\lambda_q} \right]
\quad \text{and} \\
J_m &= \left[ -\frac{2}{(2m-1)\lambda_q}, -\frac{2}{(2m+1)\lambda_q} \right], \quad m = 2,3,\ldots,
\end{split}
\end{equation}
and set $J_{-m}:= -J_m$ for $m \in \N$.
\smallskip
For even $q$, $q \geq 4$, the points in $\orbit{-\frac{\lambda_q}{2}}$ do not fall onto a boundary point
of any of the intervals $J_m$, $m \in \N$.
Indeed from the regular $\lambda_q$-CF of $-\frac{\lambda_q}{2}$ in \eqref{B3.2} and the order ``$\prec$'' in \S\ref{B4} one sees easily that
\[
-\frac{\lambda_q}{2} = \phi_0 < \phi_1 < \ldots < \phi_{{\kappa_q}-1} < -\frac{2}{3\lambda_q} <
\phi_{{\kappa_q}} = 0
\]
with $\phi_i = f_q^i\left( - \frac{\lambda_q}{2}\right)$.
If we hence define the intervals $J_{\pm1_i}$ as
\begin{equation}
\label{C3.10}
J_{\eps 1_i} := J_{\eps 1} \cap \Phi_{\eps i}
\quad \text{for } \eps = +,- ,\quad  1 \leq i \leq {\kappa_q}
\end{equation}
and therefore $J_{\eps 1_i}= \Phi_{\eps i}$ for $1 \leq i \leq {\kappa_q}-1$
we get  the partition $\cM (f_q)$, defined as
\begin{equation}
\label{C3.11}
\Iq = \bigcup_{\eps = +,-}\left( \bigcup_{ i=1}^{{\kappa_q}} J_{\eps 1_i} \, \cup \,
\bigcup_{m=2}^\infty J_{\eps m}\right),
\end{equation}
which is obviously again Markovian, since
\begin{align*}
f_q \big( J_{\eps 1_{ i}} \big)
 &= J_{\eps 1_{i+1}},\quad \eps = +,-,\, i=1,\ldots,{\kappa_q}-2 , \\
f_q\big( J_{\eps1_{\kappa_q -1}}\big)
 &= J_{\eps 1_{{\kappa_q}}}\cup \quad \bigcup_{m=2}^{\infty} J_{\eps m  }, \quad \eps = +,-
\quad \text{and} \\
f_q \big( J_{\eps {{\kappa_q}}} \big)
 &= \eps \left[ 0,\frac{\lambda_q}{2} \right],\quad \eps =+,-
\quad \text{respectively} \\
f_q \big(J_{\eps m}\big)
 &= \Iq,\quad \eps = +,-,\; m = 2,3,\ldots.
\end{align*}
The maps $f_q\big|_{J_m}$ are monotone increasing with $f_q\big|_{J_m} (x) = -\frac{1}{x}-m \lambda_q$
and \linebreak $\left(f_q\big|_{J_m}\right)^{-1} (y) = -\frac{1}{y+ m \lambda_q}$ for $m = \pm 1, \pm 2,\pm 3,\ldots$.

\smallskip

Consider next the case $q$ odd, $q \geq 5$.
In this case one has, using again the regular $\lambda_q$-CF of $-\frac{\lambda_q}{2}$ in \eqref{B3.2} and the order ``$\prec$'' in \S\ref{B4},
\[
-\frac{\lambda_q}{2} = \phi_0 < \phi_1 < \ldots < \phi_{{\kappa_q}-2} < -\frac{2}{3\lambda_q} < \phi_{{\kappa_q}-1} < -\frac{2}{5\lambda_q} < \phi_{{\kappa_q}} = 0,
\]
with ${\kappa_q} =2 h_q +1$ and the $\phi_i$'s given in \eqref{C3.5}, \eqref{C3.6}.
Hence for $\eps = +, -$ one finds $\phi_{\eps i} \in J_{\eps 1}$ for $1 \leq i \leq {\kappa_q}-2$ and $\phi_{\eps
({\kappa_q}-1)} \in J_{\eps 2}$.
If we then define for $\eps = +, -$ the intervals
\begin{align}
\label{C3.12}
J_{\eps 1_i} &:= J_{\eps 1} \cap \Phi_{\eps i}
\quad   1 \leq i \leq {\kappa_q} -1\,
\text{and hence}\, J_{\eps 1_i} =\Phi_{\eps i},\,
 1 \leq i \leq {\kappa_q} -2 \\
\nonumber
J_{\eps 2, i} &:= J_{\eps 2} \cap \Phi_{\eps i}, \, i=\kappa_q-1,
\kappa_q,
\end{align}
we find that the partition $\cM (f_q)$ defined by
\begin{equation}
\label{C3.13}
\Iq = \bigcup_{\eps =+,-}\left(\bigcup_{ i=1}^{{\kappa_q}-1} J_{\eps 1_i} \, \cup \bigcup_{
i={\kappa_q}-1}^{{\kappa_q}} J_{\eps 2_i} \, \cup \, \bigcup_{ m=3}^\infty J_{\eps m} \right)
\end{equation}
is Markovian.
Indeed for $\eps = +, -$ one finds
\begin{align*}
f_q\big( \Phi_{\eps 2i}\big)
&= \Phi_{\eps (2i+2)},\, 1 \leq i \leq h_q-2,\\
f_q \left( J_{\eps 1_{{\kappa_q} -1}}\right)
&= \eps \left[0,\frac{\lambda_q}{2}\right],\\
f_q \left( J_{\eps 2_{\kappa_q-1}} \right)
&= \eps \left[ -\frac{\lambda_q}{2},\phi_1 \right],\\
f_q \left( \Phi_{\eps (2i-1)} \right)
&= \Phi_{\eps(2i+1)} \quad \text{for $1 \leq i \leq h_q$},\\
f_q \left(J_{\eps 2_{\kappa_q}} \right) &=
\eps \left[ \phi_1,\,\, \frac{\lambda_q}{2} \right] \quad \text{and}\\
f_q \left( J_{\eps m} \right) &= \Iq  \quad\text{for $m = 3,4,\ldots$}.
\end{align*}
The maps $f_q\big|_{J_m}$ are monotone with $f_q\big|_{J_m} (x) = -\frac{1}{x}-m \lambda_q$ and $\left(f_q\big|_{J_m}\right)^{-1} (y) = -\frac{1}{y- m \lambda_q}$ for $m = \pm 1, \pm 2,\pm 3,\ldots$.

\medskip

Consider next the map $f_q^\star$ in \eqref{C1.6}.
In the case $q=3$ and $\lambda_3=1$ define the intervals $J_{\pm m}^\star,\, m=2,3,\ldots$ as
\begin{equation}
\label{C3.14}
\begin{split}
J_m^\star:=& \left[ \frac{-1}{r_3+m}, \frac{-1}{r_3+m+1} \right]
\quad \text{respectively} \\
J_{-m}^\star:=& -J_m^\star = \left[ \frac{1}{r_3+m+1}, \frac{1}{r_3+m} \right].
\end{split}
\end{equation}
Since $-R_3=-\frac{1}{2+r_3} = -1+r_3$ and $r_3= -\frac{1}{3+r_3}$, see \S\ref{B3}, we find
\begin{equation}
\label{C3.14a}
I_{R_3} = [-R_3,R_3] = \bigcup_{\eps = +, -}\quad \bigcup_{m=2}^\infty J_{\eps m}^\star
\quad \text{with} \quad
{J_m^\star}^\circ \cap {J_n^\star}^\circ = \emptyset \text{ for all }m \not=n.
\end{equation}
An easy calculation shows that
\[
f_3^\star\left( J_{\eps m}^\star \right) = \eps [r_3,R_3] \quad \text{for all } m \geq  2,
\]
where one uses $f_3^\star(-R_3) = r_3$, $\lim_{\eps \searrow 0} f_3^\star(r_3+\eps) = r_3$ and $\lim_{\eps \searrow 0} f_3^\star(r_3-\eps) = R_3$.
Hence the intervals $\left\{ J_m^\star \right\}$ define a Markov partition $\cM (f_3^\star)$ and $f_3^\star$ is a locally expanding, smooth Markov map.

\smallskip

For even $q$ define the intervals $J_{\pm m}^\star$ for $m \in \N$ as
\begin{equation}
\label{C3.15}
\begin{split}
J_m^\star:=& \left[ \frac{-1}{r_q+m\lambda_q}, \frac{-1}{r_q+(m+1)\lambda_q} \right]
\quad \text{respectively} \\
J_{-m}^\star:=& -J_m^\star = \left[ \frac{1}{r_q+(m+1)\lambda_q}, \frac{1}{r_q+m\lambda_q} \right].
\end{split}
\end{equation}
Since according to \eqref{B3.5} $R_q=1$ for $q$ even, a simple calculation shows that $f_q^\star(-R_q) = r_q = 1-\lambda_q$.
But according to \eqref{B3.8} $r_q= \lb 0 ; \ov{(1)^{h-1},2} \rb $ and hence $\left(f_q^\star\right)^{h_q-1}(r_q) = \lb 0 ; \ov{2,(1)^{h_q-1}} \rb$.
This with \eqref{B3.4} shows that
\[
 \left(f_q^\star\right)^{h_q} (R_q) = \left( f_q^\star \right)^{h_q-1}(r_q) = -\frac{1}{2\lambda_q+r_q}.
\]
The order of the points $\psi_i \in \sorbit{-R_q}$ in \eqref{C3.2} is given by
\[
-R_q = \psi_0  < \psi_1 < \ldots < \psi_{\kappa_q} = \frac{-1}{2\lambda_q + r_q}
\]
where ${\kappa_q}=h_q$ by \eqref{B3.1}.
Hence $\psi_i \in J_1^\star$ for $0 \leq i \leq {\kappa_q}-1$ whereas $\psi_{{\kappa_q}}$ is just the common boundary of $J_1^\star$ and $J_2^\star$.
Define therefore the intervals $J_{\eps 1_i}^\star$  as
\begin{equation}
\label{C3.16}
J_{\eps 1_i}^\star := J_{\eps 1}^\star \cap \Psi_{\eps i}
\quad \text{for all } 1 \leq i \leq {\kappa_q} \text{ and } \eps = +, - ,
\end{equation}
such that $J_{\eps 1}^\star = \bigcup_{i=1}^{\kappa_q} J_{\eps 1_ i}^\star$.

Then the partition $\cM (f_q^\star)$, defined by
\begin{equation}
\label{C3.17}
\Ir =\bigcup_{\eps = +,-} \left(\bigcup_{i=1}^{\kappa_q} J_{\eps 1_i}^\star \cup   \bigcup_{m=2}^\infty J_{\eps m}^\star\right),
\end{equation}
is a Markov partition, since
\begin{align*}
f_q^\star\left( J_{\eps 1_i}^\star \right)
&= J_{\eps 1_ {i+1}}^\star \quad \text{for } 1 \leq i \leq {\kappa_q} -1, \\
f_q^\star\left( J_{\eps 1_{\kappa_q}}^\star \right)
&=  \eps \left[ \frac{-1}{2\lambda_q + r_q}, R_q \right]
\quad \text{and} \\
f_q^\star\left( J_{\eps m}^\star \right)
&= \eps \left[ r_q,R_q \right]  \quad \text{for } m \geq 2.
\end{align*}
The restriction $f_q^\star \big|_{J_m^\star}$ of $f_q^\star$ to the interval $J_m^\star$ is given by
\[
f_q^\star\big|_{J_m} (x) = \frac{-1}{x} - m\lambda_q
\qquad \text{for } m \in \Z_{\not=0}
\]
and its inverse by
\[
\left(f_q^\star\big|_{J_m}\right)^{-1} (y) = \frac{-1}{y+ m\lambda_q} \qquad \text{for } y \in f_q^\star(J_m).
\]

\smallskip

Remains the case  $q$ odd, $q \geq 5$.
The intervals $J_{\eps m}^\star$ for $\eps = +, -$, $m \geq 2$ are defined as for $q$ even in \eqref{C3.15}:
\begin{equation}
\label{C3.18}
J_m^\star:= \left[ \frac{-1}{r_q+m\lambda_q}, \frac{-1}{r_q+(m+1)\lambda_q} \right]
\quad \text{and} \quad
J_{-m}^\star:= -J_m^\star.
\end{equation}
The intervals $J_{\eps 1}^\star$ are defined as
\begin{equation}
\label{C3.19}
\begin{split}
J_1^\star &:= \left[ -R_q, \frac{-1}{r_q+2\lambda_q} \right] \quad \text{respectively} \\
J_{-1}^\star &:= - J_1^\star = \left[ \frac{1}{r_q+2\lambda_q} , R_q \right].
\end{split}
\end{equation}
According to \eqref{C3.5}
\[
\psi_{2{h_q}} = \left(f_q^\star\right)^{h_q} (-R_q)  = \frac{-1}{r_q+2\lambda_q} = \lbd 0; \ov{2,1^{h_q},2,1^{h_q-1}} \rbd
\]
and
\begin{align*}
\psi_{2h_q+1}
&=  f_q^\star (-R_q) = \lbd 0 ; 2,\ov{1^{h_q-1},2,1^{h_q},2} \rbd \\
&\leq \lbd 0;3, \ov{1^{h_q},2,1^{h_q-1},2} \rbd = \frac{-1}{3\lambda_q + r_q}.
\end{align*}
Hence $\psi_i \in J_1^\star = \left[ -R_q, \frac{-1}{r_q+2\lambda_q} \right]$ for $1 \leq i \leq 2h_q= {\kappa_q}-1$ whereas $\psi_{2h_q+1} = \psi_{{\kappa_q}} \in J_2^\star$.
Define for $\eps = +, -$ the intervals $J_{\eps 1_i}^\star$  as
\begin{equation}
\label{C3.20}
J_{\eps 1_i}^\star := J_{\eps 1}^\star \cap  \Psi_{\eps i} = \Psi_{\eps i}
\quad \text{for } 1 \leq i \leq {\kappa_q}-1
\end{equation}
and the intervals $J_{\eps 2_i}^\star$ for $i=\kappa_q , \kappa_q+1$ as
\begin{align*}\label{C3.21}
J_{\eps 2_{\kappa_q}}^\star &:= J_{\eps2 }^\star \cap \Psi_{\eps \kappa_q}  = \Psi_{\eps {\kappa_q}} \quad \text{respectively}\\
\nonumber
J_{\eps 2_{\kappa_q+1}}^\star &:= J_{\eps2 }^\star \smallsetminus \Psi_{\eps \kappa_q} = \eps \left[ \psi_{\kappa_q}, \frac{-1}{3\lambda_q +r_q} \right].
\end{align*}
Then the partition
\begin{equation}
\label{C3.22}
\Ir =\bigcup_{\eps = +,-}\left(\bigcup_{i=1}^{{\kappa_q}-1} J_{\eps 1_ i}^\star \cup
\bigcup_{i=\kappa_q}^{\kappa_q+1} J_{\eps 2_ i}^\star \cup \bigcup_{m=3}^\infty J_{\eps m}^\star\right)
\end{equation}
is a Markov partition.
This follows from the following identities, which can be easily verified:
\[
f_q^\star\left(J_{\eps m}^\star \right) = \eps \left[r_q,R_q \right]
\quad \text{for all } m = 3,4,\ldots,
\]
\[
f_q^\star\left(J_{\eps 1_{2i}}^\star \right) = J_{\eps 1_{2i+2}}^\star
\quad \text{for all} 1 \leq i \leq h_q-1,
\]
\[
f_q^\star\left(J_{\eps 1_{2i-1}}^\star \right) = J_{\eps 1_{2i+1}}^\star
\quad \text{for all }1 \leq i \leq h_q,
\]
\[
f_q^\star\left(J_{\eps 1_{2h}}^\star \right)
=
\eps \left[ \psi_{\kappa_q},R_q \right] =
J_{\eps 2_{{\kappa_q}+1}}^\star \cup \bigcup_{\delta =+,-}\bigcup_{m=3}^\infty J_{\delta m}^\star \cup  J_{-\eps 2m}^\star ,
\]
\[
f_q^\star\left(J_{\eps 2_{\kappa_q}}^\star \right) =  J_{\eps 1_2}
 \quad \text{and}
\]
\[
f_q^\star\left(J_{\eps 2_{{\kappa_q}+1}}^\star \right)
=
\eps \left[ \psi_2,R_q \right]
=
\bigcup_{i=3}^{{\kappa_q}-1} J_{\eps 1_i} \cup  \bigcup_{i=\kappa_q}^{{\kappa_q}+1}J_{\eps 2_i}
\cup \bigcup_{\delta =+,-} \bigcup_{ m=3}^\infty J_{\delta m}^\star \cup J_{-\eps 2} \cup J_{-\eps 1}.
\]

\section{The maps $f_q$ and $f_q^\star$ and regular respectively dual regular $\lambda_q$-CF's}
\label{D}
We are going to use the Markov partitions $\cM(f_q)$ respectively $\cM(f_q^\star)$ constructed in the forgoing section for the maps $f_q \colon \Iq\to \Iq$ and $f_q^\star: \Ir \to \Ir$ to show that these maps can be conjugated to subshifts over infinite alphabets. By introducing sofic systems closely related to these subshifts the symbolic dynamics of the above two maps are directly related to the regular respectively dual regular $\lambda_q$-CF's.

\subsection{Symbolic dynamics for $f_q$ and a subshift of infinite type}
\label{D2}

\smallskip

For $q=3$ and $f_3 \colon I_3\to I_3$ let $F$ be the alphabet  $F=\Z \smallsetminus\{0,\pm 1\}$.
Define the transition matrix $\A=\left(\A_{i,j}\right)_{i,j\in F}$ with $\A_{i,j}\in\{0,1\}$ for $\eps = +,-$ as follows:
\begin{equation}
\label{D2.1}
\begin{split}
\A_{\eps 2,\eps m}   & = 0, \quad m\geq 2,\\
\A_{\eps 2,-\eps m}  & = 1, \quad m\geq 2,\\
\A_{\eps k, m}       & = 1, \quad k\geq 3 \quad \text{and all }  m\in F.
\end{split}
\end{equation}
Denote by $(F_{\A}^{\N}, \tau)$ the subshift over the alphabet $F$ with
\[
F_{\A}^{\N}
=
\left\{ \underline{\xi}=(\xi_i)_{i\in \N},\,\xi_i\in F, \A_{\xi_i,\xi_{i+1}}=1,\,i\in\N \right\}
\]
and $\big(\tau(\underline{\xi})\big)_i=\xi_{i+1}$ the shift map.

Let $\partial\cM(f_3):= \left\{x\in I_q: \, \exists \, n \in \Z_{\geq 0}: f_q^n(x)=0 \right\}$.
The projection map $\pi \colon I_3 \smallsetminus \partial\cM(f_3)\to F_{\A}^{\N}$ defined by
\[
\pi(x) = \underline{\xi}=(\xi_i)_{i\in \N} \quad \text{if } f_3^{i-1}(x)\in J_{\xi_i}
\quad \text{for }  i\in \N
\]
is bijective with inverse $\pi^{-1}(\underline{\xi})=x$, where $x$ is the unique point with $x \in J_{\xi_1}\cap \bigcap_{l=1}^\infty (f_{\xi_l}\circ\ldots \circ f_{\xi_1})^{-1}J_{\xi_{l+1}}$.
That the point $x$ is uniquely defined follows from the expanding property of the local branches $f_{J_m}=f_3 \big|_{J_m}$ of the map $f_3$, given on the interval $J_m$ by $f_{J_m}= -\frac{1}{x}-m$, $m \in \Z\smallsetminus\{0,\pm 1\}$.
Obviously one has $\pi \circ f_3= \tau \circ \pi$ on $I_3 \smallsetminus \partial\cM(f_3)$.

\smallskip

For $q$ even with $q=2h_q + 2$ define the alphabet $F$ as $F=\{\eps 1_i,\, \eps=+,-,\; 1\leq i\leq \kappa_q\} \cup \Z\smallsetminus\{0,\pm 1\}$.
The transition matrix $\A=\left(\A_{i,j}\right)_{i,j\in F}$ is defined in this case as follows:
\begin{equation}
\label{D2.2}
\begin{split}
\A_{\eps 1_l,\eps 1_{l+1}}       &= 1, \quad 1\leq l\leq \kappa_q-1,\\
\A_{\eps 1_{\kappa_q-1},\eps m}  &= 1, \quad m=2,3,\ldots,\\
\A_{\eps 1_{\kappa_q},-\eps 1_l} &= 1, \quad 1\leq l \leq \kappa_q,\\
\A_{\eps 1_{\kappa_q},-m}        &= 1, \quad m=2,3,\ldots ,\\
\A_{m,n}                         &= 1, \quad m\in \Z\smallsetminus\{0,\pm 1\},\quad n\in F,
\end{split}
\end{equation}
and all the other matrix elements vanishing.

Define the set $\partial\cM(f_q)$ and the map $\pi \colon \Iq \smallsetminus \partial\cM(f_q) \to F_{\A}^{\N}$ in analogy to the case $q=3$.
The same arguments as there show that this map is bijective and conjugates $f_q$ to the shift map $\tau$ with $\pi\circ f_3= \tau \circ \pi$ on $\Iq \smallsetminus \partial\cM(f_q)$.

\smallskip

For $q=2h_q+3=\kappa_q+1$ finally define the alphabet $F$ as $F=\{\eps 1_i,\, \eps=+,-,\, 1\leq i\leq \kappa_q-1\}\cup \{\eps 2_i,\, \eps=+,-,\, \kappa_q-1\leq i\leq \kappa_q\} \cup \Z\smallsetminus\{0,\pm 1,\pm 2\} $.
The transition matrix $\A=\left(\A_{i,j}\right)_{i,j\in F}$ is given in Table~\ref{D2.3}.
\begin{table}
\begin{align*}
\A_{\eps 1_{2i},\eps 1_{2i+1}}           &= 1\quad 1\leq i\leq h_q-2,\\
\A_{\eps 1_{2h_q-2},\eps 1_{2h_q}}       &= 1, \\
\A_{\eps 1_{2h_q-2},\eps 2_{\kappa_q}}   &= 1, \\
\A_{\eps 1_{2h_q},-\eps 1_{i}}           &= 1 \quad 1\leq i\leq \kappa_q-1,\\
\A_{\eps 1_{2h_q},-\eps 2_{i}}           &= 1 \quad \kappa_q\leq i\leq \kappa_q+1,\\
\A_{\eps 1_{2h_q},-\eps m}               &= 1 \quad m\geq 3,\\
\A_{\eps 1_{2i-1},\eps 1_{2i+1}}         &= 1,\quad 1\leq i\leq h_q-1,\\
\A_{\eps 1_{2h_q-1},\eps 2_{\kappa_q+1}} &= 1,\\
\A_{\eps 1_{2h_q-1},-\eps m}             &= 1,\quad m\geq 3, \\
\A_{\eps 2_{\kappa_q},-\eps 1_1}         &= 1,\\
\A_{\eps 2_{\kappa_q+1},\delta 1_i}      &= 1,\quad \delta = +,-, \quad 2\leq i \leq \kappa_q-1 ,\\
\A_{\eps 2_{\kappa_q+1},\delta 2_i}      &= 1,\quad \delta = +,- ,\quad 2\leq i \leq \kappa_q-1 ,\\
\A_{\eps 2_{\kappa_q+1},\delta n}        &= 1,\quad \delta = +,- ,\quad n\geq 3,\\
\A_{\eps 2_{\kappa_q+1},-\eps 1_1}       &= 1, \\
\A_{m,n}                                 &= 1, \quad m\in\Z\smallsetminus\{0,\pm 1,\pm 2\} ,\quad n\in F,
\end{align*}
\caption{The transition matrix $\A=(\A_{i,j})_{i,j\in F}$ with $\eps = +,-$ for $q$-regular sequences and $q$ odd, $q \geq 5$. All other matrix elements vanish.}
\label{D2.3}
\end{table}

The set $\partial\cM(f_q)$ and the map $\pi \colon \Iq \smallsetminus \partial\cM(f_q)\to F_{\A}^{\N}$  are defined similarly as in the foregoing cases $q=3$ and $q$ even and have the same properties.
The inverse $(\pi)^{-1} \colon F_{\A}^{\N}\to  \Iq \smallsetminus \partial\cM(f_q)$ is given by $\pi^{-1}(\underline{\xi})=x$ with $x \in J_{\xi_1}\cap \bigcap_{l=1}^\infty (f_{\xi_l}\circ\ldots
\circ f_{\xi_1})^{-1} J_{\xi_{l+1}}$ where $f_{\xi_i}=f\big|_{J_{\eps 1}}$ for $\xi_i=\eps 1_l$, $1 \leq l \leq \kappa_q-1$ respectively $f_{\xi_i}=f\big|_{J_{\eps 2}}$ for $\xi_i=\eps 2_l$, $\kappa_q\leq l\leq
\kappa_q+1$.
Hence also in this case the map $f_q$ gets conjugated by $\pi$ to the shift map $\tau$ on the space $ F_{\A}^{\N}$ of symbol sequences and therefore is itself a subshift of infinite type.

\subsection{A sofic system related to the map $f_q$ and the regular $\lambda_q$-CF}
\label{D2a}
The transition matrix $\A$ in \eqref{D2.1} for the subshift $f_3:I_3\to I_3$ shows that a symbol sequence $\underline{a}=(a_i)_{i\in \N}\in F_{\A}^{\N} = \pi\big(I_3\smallsetminus\partial\cM(f_3)\big)$ if and only if $(a_i,a_{i+1}) \neq(\eps 2, \eps m),\, m\geq 2$ for all $i\in \N$.
Hence this sequence is $q$-regular for $q=3$ and $F_{\A}^{\N}=\cA_3^{reg}$.
The inverse map $\pi^{-1} \colon F_{\A}^{\N}\to  \Iq \smallsetminus \partial\cM(f_q)$ therefore has the form $\pi^{-1}(\underline{a}) = \lb 0;a_1,a_2,\ldots\rb$.
This follows from

\begin{lemma}
\label{convergence3}
For $\underline{a}=(a_i)_{i\in \N} \in \cA_3^\text{reg}$ a $3$-regular sequence the limit \linebreak $\lim_{n\to\infty} \lb 0;a_1,a_2,\ldots,a_n\rb$ exists and defines a point $x\in \R$.
\end{lemma}

\begin{proof}
Set $x_n:=\lb 0;a_1,a_2,\ldots,a_n\rb$ and denote by $J_{(a_1,\ldots,a_n)}$ the closed interval $J_{(a_1,\ldots,a_n)} := J_{a_1} \cap \bigcap_{l=1}^{n-1}\left(f_{a_l}\circ\ldots f_{a_1} \right)^{-1} J_{a_{l+1}}$ with $f_{a_i} := f_q\big|_{J_{a_i}}$.
Obviously $x_n \in J_{(a_1,\ldots,a_n)}$.
All these intervals are nonempty and $J_{(a_1,\ldots,a_{n+1})}\subset J_{(a_1,\ldots,a_{n})}$ for all $n$.
Hence $\bigcap_{n=1}^\infty J_{(a_1,\ldots,a_n)}$ is not empty.
Because the map $f_3 $ is strictly expanding this set contains exactly one point $x$.
But this shows that $\lim_{n\to\infty} \lb 0;a_1,a_2,\ldots,a_n\rb = x$.
\end{proof}

In the case $q\neq 3$ the relation between the symbolic dynamics with respect to the Markov partitions $\cM(f_q)$ and the $\lambda_q$-CF is more complicated.
Indeed one has to introduce a corresponding sofic system, namely in the alphabet $F$ the letters $\eps 1_i$ respectively $\eps 2_i$ have to be replaced by the letters $\eps 1$ respectively $\eps 2$ for all
$i$.
This corresponds to replacing the Markov partition $\cM(f_q)$ defined in Section~\ref{C3} by the partition $\cJ(f_q)$ defined as $\Iq=\bigcup_{\eps=+,-}\bigcup_{m=1}^\infty J_{\eps m}$ with $J_m$ given in \eqref{C3.9}.
It is not difficult to see that this partition is generating that means $\bigcap _{i=1}^\infty f_q^{-(i-1)} J_{m_i}$ is either empty or consists of exactly one point.
This follows again from the fact that all branches of $f_q$ are expanding.
Denote by $\partial \cJ(f_q)$ the boundary points of the intervals $J_m$ including the point $x=0$ together with all their preimages under the map $f_q$.
Since $x=0$ belongs to the orbit of $-\frac{\lambda_q}{2}$ the boundaries $\partial\cM(f_q)$ and $\partial \cJ(f_q)$ coincide.
Denote by $\hat{\pi} \colon \Iq \smallsetminus \partial \cJ(f_q) \to \hat{F}^\N$ with
$\hat{F}=\Z\smallsetminus\{0\}$ the map $\hat{\pi}(x)=\underline{a}=(a_i)_{i\in \N}$ when
$f_q^{i-1}(x)\in J_{a_i}$.
If $\pi(x)=\underline{\xi}$, then obviously $a_i=m$ if $\xi_i=m\in\Z$ and $a_i=\eps 1$ if $\xi_i=\eps 1_l $ for some $l$ respectively $a_i=\eps 2$ if $\xi_i=\eps 2_l $ for some $l$.
The following Lemma then holds

\begin{lemma}
\label{D2a.1}
The map $\hat{\pi} \colon \Iq \smallsetminus\partial \cJ(f_q) \to \Ar\subset \hat{F}^\N$ is bijective.
The inverse map $\hat{\pi}^{-1}\colon \Ar\to\Iq \smallsetminus\partial \cJ(f_q)$ is given by $\hat{\pi}^{-1}(\underline{a})=\lb 0;a_1,a_2,\ldots \rb$.
\end{lemma}

\begin{proof}
Consider first $q$ even:
From the definition of the transition matrix $\A$ in \eqref{D2.2} for the case $q=2h_q+2$ it follows that there cannot be more than $h_q$ consecutive symbols $\eps 1$ in $\underline{a}=(a_i)_{i\in \N}=\hat{\pi}(x)$ since
$(\A)_{\eps 1_{\kappa_q},\eps m}=0$ for all $m=2,3,\ldots$:
indeed $h_q$ consecutive symbols $\eps 1$ are only possible for points $x$ with $f_q^{i-1}(x)\in J_{\eps 1_i}$ for $i=1,\ldots,\kappa_q$ and $f_q^{\kappa_q}\in J_{-\eps 1_l}$ for some $1\leq l \leq \kappa_q$ or $f_q^{\kappa_q}\in J_{-\eps m}$ for some $m\geq 2$.
This shows that $\underline{a}=(a_i)_{i\in \N}=\hat{\pi}(x)$ defines a $q$-regular sequence in $\Ar$.

Given on the other hand such a $q$-regular sequence $\underline{a}=(a_i)_{i\in \N}$ there exists a unique point $x\in \Iq \smallsetminus \partial \cM(f_q)$ with $\hat{\pi}(x)=\underline{a}$:
indeed if for some $l\geq 1$ and $k\geq 0$ one has $a_l=a_{l+1}=\ldots=a_{l+k}=\eps 1$ and $a_{l+k+1}=m\neq \eps 1$ consider the sequence $\underline{\xi}\in F_{\A}^{\N} $ with $\xi_{l+k}=\eps 1_{\kappa_q},\,\xi_{l+k-1}=\eps 1_{\kappa_q-1},\ldots,\, \xi_{l}=\eps 1_{\kappa_q-l}$ if $\sign m\neq \eps$ respectively $\xi_{l+k}=\eps 1_{\kappa_q-1}$, $\xi_{l+k-1}=\eps 1_{\kappa_q-2}$, $\ldots$, $\xi_{l}=\eps 1_{\kappa_q-l-1}$ if $\sign m = \eps$, whereas $\xi_i=a_i$ for all $a_i\neq \eps 1$.
Since $k\leq \kappa_q-1$ respectively $k\leq \kappa_q-2$ in the second case, the sequence $\underline{\xi}$ belongs indeed to $F_{\A}^{\N}$ and hence there exists a point $x\in \Iq\smallsetminus \partial\cM(f_q)$ with $\pi(x)=\underline{\xi}$ and hence also $\hat{\pi}(x)=\underline{a}$.
The inverse map $\hat{\pi}^{-1}$ is again given by $\hat{\pi}^{-1}(\underline{a})=\lb 0;a_1,a_2,\ldots \rb$.
Since for any $\underline{a}\in \Ar$ there exists an unique $\underline{\xi}\in F_{\A}^{\N} $ which is related to $\underline{a}$ when replacing the symbols $\pm 1_i$ by the symbol $\pm 1$, there exists therefore  $x\in I_q\smallsetminus\partial \cJ(f_q)$ with $\pi(x)=\underline{\xi}$.
But $x_n:=\lb 0;a_1,a_2,\ldots a_n\rb\in J_{(\xi_1,\ldots,\xi_n)}$ and hence $\lim_{n\to\infty}x_n = x\in\bigcap_{n=1}^\infty J_{(\xi_1,\ldots,\xi_n)}$ and hence $\hat{\pi}^{-1}(\underline{a}) = \lb 0;a_1,a_2,\ldots \rb$.

\smallskip

The same reasoning can be applied in the case $q=2h_q+3$ odd to show that the map $\hat{\pi} \colon \Iq \smallsetminus \partial \cJ(f_q)\to \Ar\in \hat{F}^\N$ is bijective with inverse  $\hat{\pi}^{-1}(\underline{a})=\lb 0;a_1,a_2,\ldots \rb$.
\end{proof}

\subsection{Symbolic dynamics for $f_q^\star$ and a subshift of infinite type}
\label{D3}
Let us start again with the case $q=3$ and recall the Markov partition $\cM(f_3^\star)$ defined in \eqref{C3.14a} by $I_{R_3}=\bigcup_{\eps = +,-}\bigcup _{m=2}^\infty J_{\eps m}^\star$  with $J_{\eps m}^\star=\eps \left[-\frac{1}{r_3+m},-\frac{1}{r_3+m+1} \right]$.
Denote by $F$ the alphabet $F= \Z\smallsetminus \{0,\pm 1\}$ and by $\A=(\A_{i,j})_{i,j\in F}$ the transition matrix with
\begin{equation}
\label{D3.1}
(\A)_{m,n} = 1
\quad \text{for all } m,n\in F \text{ with } n \neq 2 \,\sign{m},
\end{equation}
and all the other matrix elements vanishing.
Denote by $\partial\cM(f_3^\star)$ the set
\[
\partial\cM(f_3^\star)
=
\left\{ y\in I_{R_3} :\,\exists \, n\in \Z_{\geq 0} \text{ with } \left(f_3^\star\right)^n(y) = \pm r_3
   \text{ or } \left(f_3^\star\right)^n(y) = 0 \right\}
\]
and by $\pm \underline{r}_3$ the sequence $\pm \underline{r}_3=(\pm \overline{3})$.
Then one has for $f_{\xi_i}^\star := f_3^\star \big|_{J^\star_{\xi_i}}$

\begin{lemma}
\label{f_3*subshift}
The map
\[
\pi \colon I_{R_3} \smallsetminus \partial\cM(f_3^\star)
\to
F_\A^\N\smallsetminus \left\{\underline{\xi}\in F_\A^\N :\, \exists n\in \Z_{\geq 0}\,:\, \tau^n(\underline{\xi})=\pm \underline{r}_3 \right\},
\]
given by $\pi(x)=\underline{\xi}=(\xi_i)_{i\in \N}$ if $(f_3^\star)^{i-1}\in J^\star_{\xi_i}$ for $i\in \N$, is bijective, and $\pi\circ f_3^\star=\tau \circ \pi$.
Its inverse, the map $\pi^{-1}:F_\A^\N\to I_{R_3}$ can be defined on the entire set $F_\A^\N $ and is given by $\pi^{-1}(\underline{\xi})=x $ with $x$ the unique point in $I_{R_3}$ with $x\in J^\star_{\xi_1}\bigcap_{l=1}^\infty (f_{\xi_l}^\star\circ\ldots\circ f_{\xi_1}^\star)^{-1} J^\star_{\xi_{l+1}}$.
\end{lemma}

\begin{proof}
Obviously all the preimages of the point $x=0$ have a finite symbol sequence $(\xi_1,\ldots,\xi_N)$, whereas the points  $\pm r_3\in J^\star_{\pm 2} \cup J^\star_{\pm 3}$  have the two different symbol sequences $\pi (\pm r_3)=\underline{\xi} = \lbs \pm \ov{3}\rbs $ respectively $ \pi(\pm r_3) = \lbs \pm 2,\mp \ov{3}\rbs$.
The same holds then true for all the preimages of these points.
The point $x\in I_{R_3}$ is again uniquely determined because of the expansive nature of the local branches of the map $f_3^\star$.
\end{proof}

The map $f_3^\star \colon I_{R_3} \to I_{R_3}$ is hence a subshift of infinite type.

\medskip

Consider next the case $q$ even with $q = 2h_q +2$ and $\kappa_q = h_q$.
Recall the Markov partition $\cM(f_{q}^\star)$ in \eqref{C3.17}.
We define the alphabet $F$ as
\[
F=\left\{ \eps 1_i, \; \eps =+,-, \;1\leq i \leq \kappa_q \right\} \cup \Z\smallsetminus \{0,\pm 1\}
\]
and by $\A=(\A_{i,j})_{i,j\in F}$ the transition matrix with matrix elements
\begin{equation}
\label{D3.2}
\begin{split}
(\A)_{{\eps 1_i},\eps 1_{i+1}}
&= 1 \quad\text{for}\quad \eps =+,-;\quad 1\leq i \leq
\kappa_q-1,\\
(\A)_{\eps 1_{\kappa_q},m}
&= 1 \quad\text{for}\quad m\in F,\, m \neq \eps 1_i\quad 1\leq i
\leq \kappa_q, \quad \text{and}\\
(\A)_{m,n}
&= 1 \quad\text{for}\quad \abs{m} \geq 2
 \quad\text{and all} \quad n\neq \sign m 1_1,
\end{split}
\end{equation}
whereas all other matrix elements vanish.
If we define again $\partial\cM(f_q^\star)$ by
\[
\partial\cM(f_q^\star)
=
\left\{ y\in \Ir \,:\,\exists \, n\in \Z_{\geq 0} \text{ with }
(f_q^\star)^n(y) = \pm r_q \text{ or }  (f_q^\star)^n(y)=0 \right\}
\]
one shows in complete analogy with Lemma~\ref{f_3*subshift}.

\begin{lemma}
\label{D3.4}
The map
\[
\pi \colon \Ir \smallsetminus \cM(f_q^\star) \to
F_\A^\N \smallsetminus \left\{ \underline{\xi}\in F_\A^\N :\,\exists \, n\in\Z_{\geq 0}:\,  \tau^n(\underline{\xi})= \pm \underline{r}_q \right\},
\]
given by $\pi(x)=\underline{\xi}=(\xi_i)_{i\in \N}$ if $(f_q^\star)^{i-1}\in J^\star_{\xi_i}$ for $i\in \N$, is bijective, and $\pi\circ f_q^\star=\tau \circ \pi$.
Its inverse, the map $\pi^{-1}:F_\A^\N\to I_{R_q}$ can be defined on the entire set $F_\A^\N $ and is given by $\pi^{-1}(\underline{\xi})=x$ with $x$ the unique point in $\Ir$ with $x\in J^\star_{\xi_1}\bigcap_{l=1}^\infty (f_{\xi_l}^\star\circ\ldots\circ f_{\xi_1}^\star)^{-1}J^\star_{\xi_{l+1}}$.
\end{lemma}

This shows that the map $f_q^\star \colon \Ir\to \Ir$ is a subshift of infinite type also for even $q$.

\medskip

Consider finally the case $q=2h_q+3$ and $\kappa_q=2h_q+1$.
The Markov partition $\cM(f_q^\star)$ was given in this case in \eqref{C3.22}.
Define the alphabet $F$ as
\begin{align*}
F = &\left\{\eps 1_i,\; \eps = +,-,\; 1\leq i \leq \kappa_q - 1 \right\} \cup \\
&\qquad \cup \left\{\eps 2_i,\; \eps =+,-,\; \kappa_q\leq i \leq \kappa_q+1 \right\} \cup \Z \smallsetminus \{0,\pm 1,\pm 2\}.
\end{align*}
The transition matrix $\A=(\A_{i,j})_{i,j\in F}$ has now the form given in Table~\ref{D3.3}.
\begin{table}   
\begin{align*}  
(\A)_{\eps 1_{2i-1},\eps 1_{2i+1}}       &= 1\quad \text{for}\quad 1\leq i\leq h_q -1,\\
(\A)_{\eps 1_{2h_q-1},\eps 2_{\kappa_q}} &= 1\\
(\A)_{\eps 1_{2i},\eps 1_{2i+2}}         &= 1\quad \text{for}\quad 1\leq i\leq h_q -1,\\
(\A)_{\eps 1_{2h_q},\eps 2_{\kappa_q+1}} &= 1\\
(\A)_{\eps 1_{2h_q},m}                   &= 1 \quad m\in \Z\smallsetminus \{0,\pm 1,\pm 2\},\\
(\A)_{\eps 1_{2h_q},-2}                  &= 1, \\
(\A)_{\eps 2_{\kappa_q},\eps 1_2}        &= 1,\\
(\A)_{\eps 2_{\kappa_q+1},\eps 1_i}      &= 1\quad \text{for}\quad 2\leq i\leq \kappa_q -1,\\
(\A)_{\eps 2_{\kappa_q+1},\eps 2_i}      &= 1\quad \text{for}\quad  i = \kappa_q,\, \kappa_q +1,\\
(\A)_{\eps 2_{\kappa_q+1},m}             &= 1 \quad \text{for}\quad m\in \Z\smallsetminus\{0,\pm 1,\pm
2\},\\
(\A)_{\eps 2_{\kappa_q+1},-\eps 1_i}     &= 1 \quad \text{for}\quad 1\leq i\leq \kappa_q -1,\\
(\A)_{\eps 2_{\kappa_q+1},-\eps 2_i}     &= 1 \quad \text{for}\quad  i=\kappa_q ,\, \kappa_q +1,\\
(\A)_{\eps m,\eps 1_i}                   &= 1 \quad \text{for}\quad 2\leq i\leq \kappa_q -1,\quad \eps = +.-,\quad
m\in \N\smallsetminus\{0,1,2\},\\
(\A)_{\eps m,-\eps 1_1}                  &= 1,\quad \eps = +,-,\\
(\A)_{\eps m,-\eps 2_i}                  &= 1,\quad \eps = +,-,\quad I=\kappa_q ,\, \kappa_q +1.
\end{align*}    

\caption{The transition matrix $\A=(\A_{i,j})_{i,j\in F}$ with $\eps = +,-$ for $q$-dual regular sequences and $q$ odd, $q \geq 5$. All other matrix elements vanish.}
\label{D3.3}
\end{table}     

If $\partial \cM(f_q^\star)$ denotes again the set of preimages of the points $x=\pm r_q$ and the point $x=0$ one shows as in the former cases that the map
\[
\pi \colon \Ir\smallsetminus \partial \cM(f_q^\star) \to
F_{\A}^{\N} \smallsetminus \left\{\underline{\xi}\in F_{\A}^{\N}:\,\exists n\in \Z_{\geq 0}:\, \tau^n ( \underline{\xi} ) = \pm \underline{r}_q \right\}
\]
is bijective and the map $f_q^\star$ is conjugated therefore on $\Ir\smallsetminus \partial \cM(f_q^\star) $ to the shift $\tau$ on $F_{\A}^{\N}\smallsetminus \left\{ \underline{\xi} \in F_{\A}^{\N}:\, \exists n\in \Z_{\geq 0} :\, \tau^n ( \underline{\xi})= \pm \underline{r}_q \right\}$.

Hence also in the case $q$ is odd the map $f_q^\star$ is conjugate to a subshift of infinite type.

\subsection{A sofic system related to $f_q^\star$ and the dual regular $\lambda_q$-CF}
\label{D3a}
In the case of the map $f_3^\star$ the subshift $\tau: \,F_\A^\N \to F_\A^\N $ can be easily related to the dual $\lambda_3$-CF:
from the form of the transition matrix $(\A)$ in \eqref{D3.1} it follows that the sequence $\underline{b}\in F_\A^\N$ with $\underline{b}=\pi(x)$ can be characterized  by the fact that $(b_i,b_{i+1})\neq (m, 2\,\sign m)$ for all $i \in \N$ and hence $\underline{b}\in \cA_q^\text{dreg}$.
On the other hand any such sequence $\underline{b}\in \cA_3^\text{dreg}$ belongs to $F_\A^\N$ and defines a unique point $x\in \Ir$ through $J^\star_{b_1}\bigcap_{l=1}^\infty (f_{b_l}^\star\circ\ldots\circ f_{b_1}^\star)^{-1}J^\star_{b_{l+1}}$.
Since $x_n:=\lbd 0;b_1,\ldots,b_n \rbd \in J^\star_{b_1}\bigcap_{l=1}^n (f_{b_l}^\star\circ\ldots\circ f_{b_1}^\star)^{-1}J^\star_{b_{l+1}}$ for all $n$ we find $\lim_{ n\to\infty} x_n = x$ and hence $\pi^{-1}(\underline{b})=\lbd 0;b_1,b_2,\ldots,\rbd$.

To connect the subshift for the map $f_q^\star$ in the case $q=2h_q+2$ with the dual $\lambda_q$-CF one has to introduce the sofic systems by replacing in $\underline{\xi}$ all the symbols $\pm 1_i,\, 1\leq i \leq \kappa_q$ by the symbol $\pm 1$.
This corresponds to replacing the Markov partition $\cM(f_q^\star)$ by the generating partition $\cJ(f_q^\star)$ determined by $\Ir=\bigcup_{\eps =+,-}\bigcup_{m=1}^\infty J^\star_{\eps m}$ with $J^\star_{\eps m}$ defined in \eqref{C3.15}.
Denote by $\partial \cJ(f_q^\star)$ the set
\[
\partial \cJ(f_q^\star)
=
\left\{y\in I_{R_q},\,:\,\exists \, n\in \Z_{\geq 0} \text{ with } (f_q^\star)^n(y)=\pm r_q \text{ or } (f_q^\star)^n(y) = 0 \right\}
\]
which obviously coincides with the set $\partial\cM(f_q^\star)$.
Then for the alphabet $\hat{F}= \Z \smallsetminus\{0\}$ one  shows again

\begin{lemma}
\label{convergenceeta*}
The map
\[
\hat{\pi} : \Ir\smallsetminus \partial \cJ(f_q^\star)\to
\Adr \smallsetminus \left\{ \underline{b}\in\Adr:\, \exists \, n\in \Z_{\geq 0} :\, \tau^n(\underline{b})= \pm \underline{r}_q \right\} \subset \hat{F}^\N
\]
defined by
$\hat{\pi}(x)=\underline{b}=(b_i)_{i\in \N}$ if $(f_q^\star)^{i-1}\in J^\star_{b_i}$ for $i\in \N$, is bijective, and $\hat{\pi}\circ f_q^\star=\tau \circ \hat{\pi}$.
Its inverse, the map $\hat{\pi}^{-1}:\Adr\to I_{R_q}$
can be defined on the entire set $\Adr $ and is given by $\hat{\pi}^{-1}:(\underline{b})=\lbd 0;b_1,b_2,\ldots,\rbd$.
\end{lemma}

\begin{proof}
Since $\pm \psi_{\kappa_q}=(f_q^\star)^{h_q-1}(\pm r_q)\in J_{\pm 1}^\star\cap J^\star_{\pm 2}$ this point has two different dual regular sequences $\underline{b}=\pm (1,(-1)^{h_q},\ov{-2,(-1)^{h_q}-1})$ respectively $\underline{b}=\pm (\ov{1^{h_q-1},2})$.
Hence also all preimages of this point have two different dual regular sequences, but these points all have the same tail as the point $\pm r_q$.
If $\underline{b}=(b_i)_{i\in\N}=\hat{\pi}(x)$ assume that $\underline{b}$ contains for some $k \geq 0$ and some $l\geq 0$ a subsequence $b_{k+1}=\dots = b_{k+l}=\eps 1$ with either $k=0$ or $b_k \neq \eps 1$ and $b_{k+l+1}\neq \eps 1$.
Then the sequence $\underline{\xi}=(\xi_i)_{i\in\N}$ related to $\underline{b}$ by replacing the symbols $\pm 1_i$ by the symbol $\pm 1$ must be of the form $\xi_{k+l}=\eps 1_{\kappa_q}$ and hence $\xi_{k+i}=\eps 1_{\kappa_q-(l-i)}$ for $1\leq i \leq l$.
This shows that $l\leq \kappa_q$.
The case $l=\kappa_q$ however is only possible if either $k=0$ or $b_k=-\eps m$, $m \geq 1$.
This shows that $\underline{b}=\hat{\pi}(x)\in \Adr$.

Given on the other hand  $\underline{b}\in \Adr$ with a subsequence $b_{k+1}=\ldots=b_{k+l}=\eps 1$ and $b_{k+l+1}\neq \eps 1$ for some $k$ and some $l$ then define the sequence $\underline{\xi}$ such that $\xi_{k+i}=\eps 1_{\kappa_q-(l-i)}$, $1\leq i\leq l$.
Since $l\leq \kappa_q $ respectively $l\leq \kappa_q-1$ the sequence $\underline{\xi}$ belongs to $F_{\A}^{\N}$ and hence there exists a point $x\in \Ir$ with $\pi(x)=\underline{\xi}$ and therefore by construction also $\hat{\pi}(x)=\underline{b}$.
The inverse map $\hat{\pi}^{-1}$ is obviously defined for all $\underline{b}\in \Adr$.
An argument completely analogous to the one in the case $q=3$ then shows that $\hat{\pi}^{-1}(\underline{b}) = \lbd 0;b_1,b_2,\ldots \rbd$.
\end{proof}

\smallskip

Introduce finally in the case $q=2h_q+3$ for the map $f_q^\star$ the sofic system defined by replacing in $\underline{\xi}\in F_\A^\N$ all the symbols $\pm 1_i$ by the symbol $\pm 1$ and the symbols $\pm 2_i$ by the symbol $\pm 2$.
Denote by $\cJ(f_q^\star)$ the corresponding generating partition $\Ir=\bigcup_{\eps =+,-}\bigcup_{m=1}^\infty J^\star_{\eps m}$ and by $\partial \cJ(f_q^\star)$ the set of preimages of the points $x=\pm r_q$ and $x=0$ which obviously coincides with the set $\partial \cM(f_q^\star)$.
As in the previous cases one shows also for $q=2h_q+3$

\begin{lemma}
The map
\[
\hat{\pi}:\Ir\smallsetminus \partial \cJ(f_q^\star)\to
\Adr\smallsetminus \left\{\underline{b}\in\Adr : \, \underline{b} \text{ has the tail } \pm \underline{r}_q \right\} \subset \hat{F}^\N
\]
defined by
$\hat{\pi}(x)=\underline{b}=(b_i)_{i\in \N}$ if $(f_q^\star)^{i-1}\in J^\star_{b_i}$ for $i\in \N$, is bijective and $\hat{\pi}\circ f_q^\star=\tau \circ \hat{\pi}$.
Its inverse, the map $\hat{\pi}^{-1} \colon \Adr\to \Ir$
can be defined on the entire set $\Adr $ and is given by $\hat{\pi}^{-1} \colon(\underline{b}) = \lbd
0;b_1,b_2,\ldots \rbd$.
\end{lemma}

\begin{proof}
From the form of the transition matrix $\A$ in Table~\ref{D3.3} it is clear that there are only restrictions on the symbol sequence $\underline{b}$ for $\underline{b}=(b_i)_{i\in \N}=\hat{\pi}(x)$ if it contains subsequences of consecutive  symbols $\pm 1$ since $(\A)_{i,j}=1$ for all $j\in F$ if $\abs{i}\geq 3$.
Assume $b_k=  m \neq 1$ and $b_{k+1}=\ldots=b_{k+l}=\pm 1,\, b_{k+l+1}\neq \pm 1$ for some $k\geq 0$ and some $l\geq 1$, where $k=0$ means that $b_1=\pm 1$.
Then either $(f_q^\star)^{k+l-1}(x) \in J^\star_{\pm 1_{2h_q-1}}$ or $(f_q^\star)^{k+l-1}(x)\in J^\star_{\pm 1_{2h_q}}$ and hence $\xi_{k+l}=\pm 1_{2h_q-1}$ or $\xi_{k+l}=\pm 1_{2h_q}$.
In the first case $\xi_{k+i}=\pm 1_{2h_q-2(l-i)-1},\, 1 \leq i \leq l$ and hence $\xi_{k+1}=\pm 1_{2h_q-(2l-1)}$.
If $m= \pm n$ for some $n\geq 3$ then necessarily $l\leq h_q-1$ since $(\A)_{m,\pm 1_1}=0$ for all $m\in F$.
If on the other hand $m=\mp n $ for some $n \geq 3$ or $k=0$ then $l\leq h_q$ with $l=h_q$ iff $\xi_{k+1}=\pm 1_1$.
In the case $(f_q^\star)(x)\in J^\star_{\pm 1_{2h_q}}$ we find $\xi_{k+i}=\pm 1_{2h_q-2(l-i)},\, 1 \leq i \leq l$ and hence $\xi_{k+1}=\pm 1_{2+2(h_q-l)}$.
This shows that also in this case $l\leq h_q$.
In the symbol sequence $\underline{b}$ there can appear therefore no subsequence of more than $h_q$ consecutive symbols $\pm 1$.

Assume next that there exists in $\underline{b}$ a subsequence of $h_q$ consecutive symbols $\pm 1$ such that $b_{k+1}=\ldots=b_{k+h_q}=\pm 1 $ and $b_{k+h_q+1}\neq=\pm 1$.
Then either $k=0$, that means $b_{k+1}=b_1$, or $b_k= \mp n$ for some $n\neq 1$.
Then $(f_q^\star)^{k+h_q-1}(x)\in J^\star_{\pm 1_{2h_q-1}}$ or $(f_q^\star)^{k+h_q-1}(x)\in J^\star_{\pm 1_{2h_q}}$ and hence $\xi_{k+h_q}=\pm 1_{2h_q-1}$ respectively $\xi_{k+h_q}=\pm 1_{2h_q}$.
The transition matrix $\A$ in Table~\ref{D3.3} then shows that in the first case $\xi_{k+h_q+1}= \pm 2_{\kappa_q}$ and in the second case
\[
\xi_{k+h_q+1}\in \left\{n\in F,\; n\neq \pm 1_i,\; 1 \leq i \leq \kappa_q-1, \;n \neq \pm 2_{\kappa_q}\right\}.
\]
If $\xi_{k+h_q+1}= \pm 2_{\kappa_q}$ then $\xi_{k+h_q+1+i}= \pm 1_{2i}$ for $1\leq i\leq h_q$ and hence
$\xi_{k+2h_q+1}= \pm 1_{2h_q}$.

If in the second case $\xi_{k+h_q+1}= \pm 2_{\kappa_q+1}$ then the maximal number of consecutive symbols $\pm 1$ in $\underline{b}$ is $h_q-1$ since in this case $\xi_{k+h_q+2}=\pm 1_i$ for some $i \geq 3$ and only for $i=3$ one has $\xi_{k+2h_q}=\pm 1_{2h_q-1}$.
In all other cases when $b_{k+h_q+1}\neq \pm 2$ the number of consecutive symbols $\pm 1$ is certainly bounded by $h_q$.
This shows that in the sequence $\underline{b}=\hat{\pi}(x)$ the subsequence $(\pm m, (\pm 1)^h_q,\pm 2, (\pm 1)^h_q,\pm 2)$ and the subsequence $(\pm 1)^{h_q+1}$ cannot appear.
Hence $\underline{b}=\hat{\pi}(x)\in \Adr$.
Since also in this case to every $\underline{b}\in\Adr$ there exist a unique $\underline{\xi}\in F_\A^\N$ which is related to $\underline{b}$ by replacing the symbols $\pm 1_i$ respectively the symbols $\pm 2_i$ by the symbols $\pm 1$ respectively $\pm 2$, the same arguments as in the previous cases apply to show, that the inverse map $\hat{\pi}^{-1}$ is given by $\hat{\pi}^{-1}(\underline{b}) = \lbd 0;b_1,b_2,\ldots\rbd$.
\end{proof}

\subsection{Symbolic dynamics and the natural extension $F_q$ of the map $f_q$}
\label{D1}
Consider the maps $f_q$ and $f_q^\star$.
Since
\begin{align*}
f_q\big( \lb 0;a_1,a_2,a_3,a_4,\ldots \rb\big) &= \lb 0;a_2,a_3,a_4,\ldots \rb
\quad \text{and} \\
f_q^\star\big( \lbd 0;a_0,a_{-1},a_{-2},a_{-3},\ldots \rbd\big) &= \lbd 0;a_{-1},a_{-2},a_{-3}, \ldots \rbd ,
\end{align*}
$f_q$ and $f_q^\star$ are equivalent to the shift map $\tau$ on one sided infinite sequences $\underline{a}_>:=(a_i)_{i\in \N}\in \Ar$ respectively $\underline{a}_< :=(a_i)_{i \in \Z_{\leq 0}}\in \Adr$.
Denote by $\cA_q$ the set of two-sided infinite sequences
\[
\cA_q =
\left\{ \underline{a}=(a_i)_{i\in\Z}: \, \forall l\in\Z ,\,\forall k>0 \,:\, (a_l,a_{l+1},\ldots,a_{l+k})\notin \cB_q \right\},
\]
where $\cB_q$ was defined in \eqref{B2.3}.
The natural extension of the one-sided shift map $\tau$ is the two-sided shift $\tau:\cA_q\to\cA_q$ with
\begin{equation}
\label{D1.8}
\begin{split}
\big(\tau(\underline{a})\big)_i&= a_{i+1},\, i\in\Z \, \quad
\text{respectively its inverse}\\
\big(\tau^{-1}(\underline{a})\big)_i&= a_{i-1},\, i\in\Z \quad \text{if}\, \underline{a}=(a_i)_{i\in\Z}.
\end{split}
\end{equation}

The natural extension $F_q$ of the map $f_q$ respectively its inverse $F_q^{-1}$ can then be identified simply with the corresponding induced maps on pairs of points $(x,y)$ with regular respectively dual regular $\lambda_q$-CF
\[
x =  \lb 0;a_1,a_2,a_3,a_4,\ldots \rb
\text{ and }
y = \lbd 0;a_0,a_{-1},a_{-2},a_{-3},a_{-4},\ldots \rbd,
\]
as long as the two-sided sequence $\underline{a}=(a_i)_{i\in\Z}$ belongs to $\cA_q$.
Then, $F_q$ and $F_q^{-1}$ satisfy
\begin{equation}
\label{D1.5}
\begin{split}
F_q \big(\lb 0;a_1,a_2,\ldots \rb, \lbd 0;a_0,\ldots \rbd \big) &= \big(\lb 0;a_2,\ldots \rb, \lbd 0;a_1,a_0,\ldots \rbd \big) \text{ and} \\
F_q^{-1} \big(\lb 0;a_1,\ldots \rb, \lbd 0;a_0,a_{-1},\ldots \rbd \big) &= \big(\lb 0;a_0,a_1,\ldots \rb, \lbd 0;a_{-1},\ldots \rbd \big).
\end{split}
\end{equation}
To characterize the set $\Omega_q$ of pairs $(x,y)$ with the above property, define in a first step $\Iqs := \Iq \smallsetminus \left\{x \text{ has a finite regular $ \lambda_q$-CF} \right\}$.
Obviously $\Iqs$ has full measure.
Denote next by $\Pi_1 \colon \cA_q \to  \Iqs$ the map
\begin{equation}
\label{D1.1}
\Pi_1(\ldots,  a_{-1}, a_0; a_1, a_2, a_3 \ldots ) = \lb 0;a_1,a_2,a_3,\ldots \rb.
\end{equation}
By construction the following lemma holds.

\begin{lemma}
\label{D1.2}
The map $\Pi_1$ is surjective and satisfies $\Pi_1 \circ f_q = \tau \circ \Pi_1$.
\end{lemma}

\smallskip

Next, define $\Irs := \Ir \smallsetminus \left\{y \text{ has a finite dual regular $\lambda_q$-CF} \right\}$ which has full measure.
Similar to \eqref{D1.1} the map $\Pi_2 \colon \cA_q \to \Irs$ given by
\begin{equation}
\label{D1.3}
\Pi_2 \lbs \ldots, a_{-2}, a_{-1}, a_0; a_1 , \ldots \rbs = \lbd 0;a_0,a_{-1},a_{-2},\ldots \rbd.
\end{equation}
is well defined, surjective and satisfies $\Pi_2 \circ f_q^\star = \tau^{-1} \circ \Pi_2$.
The following Lemma characterizes the domain of definition $\Omega_q$ of the natural extension $F_q$:
\begin{lemma}
\label{D1.4a}
For $\underline{a} \in \cA_q$ we have
\[
\begin{array}{ll}
\Pi_2(\underline{a}) \in \pm \left[ \psi_{\kappa_q-i+1},R_q \right]   \; & \mbox{if }\; \Pi_1(\underline{a}) \in \pm \Phi_i =\pm \left[\phi_{i-1},\phi_i\right], \;  i\in\{1,\ldots,\kappa_q\}.
\end{array}
\]
The set $\Omega_q\subset \Iq \times \Ir$ hence is given by
\[
\Omega_q =\bigcup_{i=1}^{\kappa_q} \Big(\big[\phi_{i-1},\phi_i\big] \times \big[\psi_{\kappa_q-i+1},R_q\big] \Big)
  \cup \Big(\big[-\phi_i,-\phi_{i-1}\big] \times \big[-R_q,-\psi_{\kappa_q-i+1}\big]  \Big).
\]
\end{lemma}

\begin{proof}
W.l.o.g.\ assume $\Pi_1(\underline{a}) \in \Phi_i \subset \Iq$ where $\Phi_i$ is defined in \eqref{C3.6}.

For $q$ even Lemma~\ref{C3.3} and the $\lambda_q$-CF of $-\frac{\lambda}{2}$ in \eqref{B3.2} show
\[
\lb 0;(1)^{h_q-i+1} \rb \preceq \Pi_1(\underline{a}) \preceq \lb 0;(1)^{h_q-i} \rb.
\]
Hence $\Pi_1(\underline{a})$ has a $\lambda_q$-CF of the form $\Pi_1(\underline{a}) = \lb 0; (1)^{h_q-i}, m, \ldots \rb$ for some $m \geq 2$.
Since $\underline{a}$ is $q$-regular we have
\[
\underline{a} = \lbs \ldots, a_{-i-1}, a_{-i}, a_{-i+1},\ldots, a_0; (1)^{h_q-i},m,a_{h_q-i+2},\ldots \rbs
\]
with at most $i-1$ consecutive $1$'s in the sequence $\lbs a_{-i+2},\ldots, a_0\rbs$.
The point $\Pi_2(\underline{a})$ hence is bounded by the largest and smallest number whose dual regular $\lambda_q$-CF starts with at most $i-1$ consecutive $1$'s and hence
\[
\lbd 0; (1)^{i-1},2,\ov{(1)^{h_q-i},2} \rbd \preceq \Pi_2(\underline{a}) \preceq \lbd 0;(-1)^{h_q},-2,\ov{(-1)^{h_q-1},-2} \rbd.
\]
But \eqref{C3.4} and \eqref{B3.4} show that these bounds are just $\psi_{\kappa_q-i+1}$ and $R_q$.

\smallskip

The case $q$ odd, $q\geq 5$, is slightly more complicated.
First, assume $i$ to be even and put $j = \frac{i}{2}$.
Then by Lemma~\ref{C3.3} and the $\lambda_q$-CF of $-\frac{\lambda}{2}$ in (\ref{B3.2})
\[
\Pi_1(\underline{a}) \in \Phi_{2j}=\Big[ \lb 0; (1)^{h_q-j+1} \rb, \lb 0;(1)^{h_q-j},2,(1)^{h_q} \rb \Big].
\]
Hence $\Pi_1(\underline{a})$ has a $\lambda_q$-CF of the form $\Pi_1(\underline{a}) = \lb 0; (1)^{h_q-j},2,(1)^{h_q}, m, \ldots \rb$ for some $m \geq 2$.
Since $\underline{a} \in \cA_q$ the sequence $(a_{-j+1},\ldots,a_0)$ in
\[
\underline{a}= \lbs \ldots, a_{-j}, a_{-j+1},\ldots, a_0; (1)^{h_q-j},2,(1)^{h_q},m ,a_{2h_q-j+2},\ldots \rbs
\]
cannot contain more than $j-1$ consecutive digits $1$.
Hence $\Pi_2(\underline{a})$ is bounded by the points
\begin{multline}
\nonumber
\lbd 0; (1)^{j-1},2,\ov{(1)^h,2,(1)^{h-1},2} \rbd \\
 \preceq \Pi_2(A) \preceq \lbd 0;(-1)^h,-2,\ov{(-1)^h,-2,(-1)^{h-1},-2} \rbd.
\end{multline}
which by Lemma~\ref{C3.3} and \eqref{B3.4} are just $\psi_{\kappa_q-2j+1}$ and $R_q$.

Next, consider the case $i$ odd and put $j = \frac{i-1}{2}$ for $1\leq i\leq h_q$.
Again, by Lemma~\ref{C3.3} and \eqref{B3.4}
\[
\Pi_1(\underline{a}) \in \Phi_{2j+1} = \Big[ \lb 0;(1)^{h_q-j},2,(1)^{h_q} \rb, \lb 0; (1)^{h_q-j} \rb \Big],
\]
and therefore $\Pi_1(\underline{a})$ has a $\lambda_q$-CF of the form $\Pi_1(\underline{a}) = \lb 0; (1)^{h_q-j}, m, \ldots \rb$ for some $m \geq 2$.
Hence there is again a restriction on the sequence $\underline{a}$:
\[
\underline{a}=\Big(\ldots, a_{-h-j-1}, \underbrace{a_{-h-j}, \ldots, a_0; (1)^{h-j},m}_{\text{not a forbidden block}},a_{h-j+1},\ldots \Big),
\]
and therefore $\Pi_2(\underline{a})$ is bounded by
\begin{multline}
\nonumber
\lbd 0; (1)^{j},2,\ov{(1)^{h_q-1},2,(1)^{h_q},2} \rbd   \preceq \Pi_2(A)  \\
 \preceq \lbd 0;(-1)^{h_q},-2,[0; (-1)^{h_q},-2,\ov{(-1)^{h_q},-2,(-1)^{h_q-1},-2} \rbd,
\end{multline}
which by Lemma~\ref{C3.3} and~\eqref{B3.2} respectively~\eqref{B3.4} are just  $\psi_{\kappa_q-2j}$ and $R_q$ .

Finally for $q=3$ Lemma~\ref{C3.3} and \eqref{B3.4} show $\lb 0;2 \rb \preceq \Pi_1(\underline{a}) \prec \lb 0; \rb,$ and therefore $\Pi_1(\underline{a})$ has a $\lambda_3$-CF of the form $\Pi_1(\underline{a}) = \lb 0; m, \ldots \rb$ for some $m \geq 2$.
Hence $\Pi_2(\underline{a})$ must not have a leading digit $2$.
This implies the bounds $r_3=[0;\ov{3}] \preceq \Pi_2(\underline{a}) \preceq [0;-2,\ov{-3}] = R_3$.

\smallskip

The case $\Pi_1(\underline{a})\in \Phi_{-i}$ for some $1\leq i \leq \kappa_q$ follows from $\Pi_1(-\underline{a}) = - \Pi_1(\underline{a})$.
\end{proof}

\medskip

Recall the definition of the domain $\Omega_q$ in Lemma~\ref{D1.4a} and define the set $\Omega_q^\star = \Omega_q \cap ( \Iqs  \times \Irs )$, which obviously is dense in $\Omega_q$.
Then one has

\begin{lemma}
\label{D1.7}
The map $\Pi \colon \cA_q \to \Omega_q^\star$ with $\Pi(\underline{a}) = \Big(\Pi_1(\underline{a}), \Pi_2(\underline{a})\Big)$ is a bijection.

For $F_q : \Omega_q^\star \to \Omega_q^\star$ and $F_q^{-1} : \Omega_q^\star \to \Omega_q^\star$ given by~\eqref{D1.5} the diagrams
\[
\begin{array}{ccc}
\cA_q & \stackrel{\tau}{\longrightarrow} & \cA_q \\
{\scriptstyle \Pi} \downarrow && \downarrow {\scriptstyle \Pi}  \\
\Omega_q^\star &\stackrel{F_q}{\longrightarrow} & \Omega_q^\star
\end{array}
\quad \text{and} \quad
\begin{array}{ccc}
\cA_q & \stackrel{\tau^{-1}}{\longrightarrow} & \cA_q\\
{\scriptstyle \Pi} \downarrow && \downarrow {\scriptstyle \Pi}  \\
\Omega_q^\star &\stackrel{F_q^{-1}}{\longrightarrow} & \Omega_q^\star
\end{array}
\quad \text{commute}.
\]
\end{lemma}

\begin{proof}
Obviously, the map $\Pi$ is well defined.
Commutativity of the diagrams follows from combining Lemma~\ref{D1.2} and
\begin{multline}
\nonumber
\Pi_2\big(\tau(\underline{a})\big) = \Pi_2\big(\tau( \ldots, a_{-1}, a_0; a_1, a_2, \ldots)\big)
=
\Pi_2\big( \ldots, a_{-1}, a_0, a_1; a_2, \ldots \big) \\
=
\lbd 0;a_1,a_0,a_{-1}, \ldots \rbd
=
\frac{-1}{\lbd a_0,a_{-1},\ldots \rbd + a_1 \lambda_q}
=
\frac{-1}{\Pi_2(\underline{a}) + a_1 \lambda_q}
\end{multline}
respectively
\begin{multline}
\nonumber
\Pi_1\big(\tau^{-1}(\underline{a})\big) = \Pi_1\big(\tau^{-1}( \ldots, a_{-1}, a_0; a_1, \ldots) \big)
=
\Pi_1\big( \ldots, a_{-1}; a_0, a_1, \ldots \big) \\
=
\lb 0;a_0,a_1,a_2, \ldots \rb
=
\frac{-1}{\lb a_1,a_2,\ldots \rb + a_0 \lambda_q}
=
\frac{-1}{\Pi_1(\underline{a}) + a_0 \lambda_q}.
\end{multline}

Since the map $\Pi \colon \cA_q \to \Omega_q^\star$ is obviously injective we only need to show $\Pi(\cA)=\Omega_q^\star$.
For this take $(x,y) \in \Omega_q^\star$ with
\[
x= \lb 0;a_1,a_2,a_3,\ldots \rb
\quad \mbox{and} \quad
y= \lbd 0;a_0,a_{-1},a_{-2},\ldots \rbd.
\]
If $(x,y) \in \Big([\phi_{i-1},\phi_i] \cap \Iqs \Big) \times \Big([\psi_{\kappa_q-i+1},\psi_0] \cap \Irs \Big)$ for some $i \in \{ 1, \ldots,\kappa_q\}$ the definitions of $\phi_i$ and $\psi_i$ as elements of $\orbit{-\frac{\lambda_q}{2}}$ and $\orbit{-R_q}$ respectively imply that the bi-infinite sequence
\[
\underline{a}:= \lbs \ldots,a_{-2},a_{-1},a_0;a_1,a_2,a_3,\ldots \rbs
\]
does not contain forbidden blocks, and hence $\underline{a} \in \cA_q$.
\end{proof}

It is well known that the map $\Pi \colon \cA_q \to \R^2$ is continuous (see for instance \cite{R04}) when $\cA_q$ is equipped with the usual metric of the shift space.

\section{Some applications}
\label{E}

\subsection{Reduced geodesics on $\Gn{q} \backslash \HH$ and the natural extension of $f_q$}
\label{E1}
The Poin\-car\'e upper half-plane, equipped with the hyperbolic metric $\dd s$ with $\dd s^2 = \frac{\dd x^2 + \dd y^2}{y^2}$, is denoted by $\HH = \{z \in \C; \; \im{z}>0\}$.
The group of isometries of this space is given by $\PSL{\R}$.
The boundary of $\HH$ is the projective line $\PP$.

We consider oriented geodesics on $\HH$.
Geodesic lines on $\HH$ are half-circles perpendicular to $\R$ or straight lines parallel to the imaginary axis $\re{z}=0$.
An oriented geodesic $\omega$ on $\HH$ will be represented by the two base points $\omega_-, \omega_+ \in \R\cup\{i\infty\}$ with its orientation from $\omega_-$ towards $\omega_+$.
We denote such a geodesic by $\omega=(\omega_-,\omega_+)$.

We call two oriented geodesics $\omega$ and $\upsilon$ \emph{$\Gn{q}$-equivalent} if there exists an element $g \in \Gn{q}$ with $g\, \omega_- = \upsilon_-$ and $g \, \omega_+ = \upsilon_+$.

Then one can show

\begin{theorem}
\label{E1.1}
Let $\omega = (\omega_-,\omega_+)$ be a geodesic with $\omega_-$ having an infinite regular respectively $\omega_+$ having an infinite dual regular $\lambda_q$-CF.
Then there exist a geodesic $\omega^\prime = (\omega^\prime_-,\omega^\prime_+)$ such that
\begin{itemize}
\item $\omega$ and $\omega^\prime$ are $\Gn{q}$-equivalent and
\item $\big(S\, \omega^\prime_+, - \omega^\prime_-\big) \in \Omega_q$.
\end{itemize}
\end{theorem}

\begin{proof}
We prove the Theorem first for  $q\geq 4$.
Using translations by powers of $T_q$ we may assume that for  $\omega=(\omega_-,\omega_+)$ either
\begin{itemize}
\item $\omega_+>0$ and $\omega_-\in \left[-R_q,-r_q\right] \subset \Ir$ or
\item $\omega_+<0$ and $\omega_-\in \left[r_q,R_q\right] \subset \Ir$.
\end{itemize}

Assume $\omega_+>0$ and $\omega_-\in\left[-R_q,-r_q\right]$ with infinite $\lambda_q$-CF's
\begin{equation}
\label{E1.2}
\omega_+= \lb  a_0;a_1,a_2,\ldots \rb
\quad \text{and} \quad
\omega_-= \lbd 0;b_1,b_2,\ldots \rbd.
\end{equation}

For $x = S \, \omega_+$ and $y= -\omega_-$ we have
\begin{align}
\label{E1.3}
x &=
\begin{cases}
[ 0;a_0,a_1,a_2,\ldots ]  & \text{if } a_0 \not=0, \\
\lb  a_1;a_2 \ldots \rb   & \text{if } a_0 =0 \text{ and}
\end{cases} \\
\nonumber
y &= \lbd 0;-b_1,-b_2,\ldots \rbd.
\end{align}
The following three cases have to be discussed: ``$a_0 \geq 2$'', ``$a_0 = 1$'' and ``$a_0 = 0$''.

If $a_0 \geq 2$ then $\omega_+ \geq  2\lambda - \frac{\lambda}{2}$ and the $\lambda_q$-CF of $x$ in \eqref{E1.3} is regular.
If the two-sided sequence
\begin{equation}
\label{E1.4}
\underline{a}:= \lbs \ldots,-b_2,-b_1;a_0,a_1,\ldots \rbs
\end{equation}
belongs to $\cA_q$ then by Lemma~\ref{D1.7} $\big(S\, \omega_+, - \omega_-\big)= (x,y) = \Pi(\underline{a}) \in \Omega_q^\star \subset \Omega_q$ .

Assume therefore $\underline{a} \not\in \cA_q$.
The forbidden block must appear around ``$;$'' in \eqref{E1.4}.
This can only happen for $q$ even if $\lbs -b_{h_q}, \ldots -b_1 \rbs = \lbs 1^{h_q} \rbs$  respectively for $q$ odd if $\lbs -b_{2h_q+1}, \ldots -b_1 \rbs = \lbs 1^{h_q},2,1^{h_q} \rbs$.
Using the lexicographic order ``$\prec$'' in Section \ref{B4} we find
\begin{multline}
\nonumber
-\omega_-
=
\lbd 0; -b_1, -b_2, \ldots \rbd
=
\begin{cases}
\lbd 0;1^{h_q},\ldots \rbd       & \text{for $q$ even and} \\
\lbd 0;1^{h_q},2,1^{h_q},\ldots \rbd & \text{for $q$ odd}
\end{cases}  \\
\prec
\begin{cases}
\lb 0; \ov{(1)^{h_q-1},2} \rb & \text{for $q$ even and} \\
\lb 0; \ov{(1)^{h_q},2,(1)^{h_q-1},2} \rb & \text{for $q$ odd, $q \geq 5$}.
\end{cases}
\;= r_q.
\end{multline}
Then Lemma~\ref{B4.4} implies $\lbd 0; -b_1, -b_2, \ldots \rbd < r_q$.

On the other hand, $-\omega_- \in \left[r_q,R_q\right]$ implies $-\omega_- \geq r_q$.
This leads to a contradiction.

\smallskip

If $\mathbf{a_0=0}$ then $\omega_+ \in  \left(0,\frac{\lambda_q}{2}\right] \subset \left(0,\frac{2}{\lambda_q}\right)$.
Hence $a_1 <0$ in the $\lambda_q$-CF \eqref{E1.2}.
For $m \in \Z_{\geq 2}$ the sequence $\lbs -m,a_1,a_2,\ldots \rbs$ is $q$-regular
and also $\lbs -m,b_1,b_2,\ldots \rbs$ is dual $q$-regular, since $r_q \leq -\omega_-$ implies
\begin{align*}
r_q &=
\begin{cases}
\lbd 0;\overline{(1)^{h_q-1},2} \rbd              & \text{if $q$ is even and} \\
\lbd 0;\overline{(1)^{h_q},2,(1)^{h_q-1},2} \rbd  & \text{if $q$ is odd.}
\end{cases} \\
&\preceq \lbd 0; -b_1, -b_2, \ldots \rbd = -\omega_-.
\end{align*}
For $g := S T_q^{-m}$ define $\big(\omega^\prime_-,\omega^\prime_+\big) = \omega^\prime := g \; \omega$ with
\begin{align*}
\omega_+^\prime &= \lb  0; -m,a_1,a_2,\ldots \rb > 0 \quad \text{and} \\
\omega^\prime_- &= \lbd 0; -m,b_1,b_2,\ldots \rbd \in (0,-r_q] \subset [r_q,R_q].
\end{align*}
The corresponding bi-infinite sequence $\underline{a}^\prime:= \lbs \ldots,-b_1,m;-m,a_1,\ldots \rbs$ is then $q$-regular and Lemma~\ref{D1.7} hence implies $\big(S\, \omega_+^\prime, - \omega_-^\prime\big) = \Pi(\underline{a}^\prime) \in \Omega_q^\star \subset \Omega_q$.

\smallskip

For $\mathbf{a_0 = 1}$
\begin{gather}
\nonumber
\omega_+ \geq  \lambda_q - \frac{\lambda_q}{2} =
\begin{cases}
\lb 1;(1)^{h_q} \rb & \text{for $q$ even,} \\
\lb 1;(1)^{h_q},2,(1)^{h_q} \rb & \text{for $q$ odd and} \\
\end{cases} \\
\label{E1.6}
-R_q \leq \omega_- \leq -r_q =
\begin{cases}
\lbd 0;\overline{(-1)^{h_q-1},-2} \rbd 		& \text{for $q$ even,} \\
\lbd 0;\overline{(-1)^{h_q},-2,(-1)^{h_q-1},-2} \rbd 	& \text{for $q$ odd.}
\end{cases}
\end{gather}
Then the dual regular $\lambda_q$-CF of $y$ cannot start for $q$ even with a block of the form $\lbs 1^{h_q} \rbs$ respectively for $q$ odd $\lbs 1^{h_q},2,1^{h_q} \rbs$.

Assume $\underline{a} \not\in \cA_q$, with $\underline{a}$  defined as in \eqref{E1.4}.
Then a block of the form $\lbs 1^l;1,1^t \rbs$ with $0 \leq l \leq h$ and $0 \leq t \leq h$ must exist around the ``$;$'' in $\underline{a}$ such that $a_{t+1} \not=1$, $-b_{l+1} \not=1$.
As in the case $a_0=0$, choose an $m \in \Z \smallsetminus\{-2,-1,0,1\}$ with $(m+1,a_1,a_2,\ldots)$ being $q$-regular and $(m,b_1,b_2,\ldots)$ being dual $q$-regular: indeed any $m$ with $\sign m =-\sign{b_1}$ can be used.
Define $g := ST_q^m$ and $\big(\omega^\prime_-,\omega^\prime_+\big) = \omega^\prime := g\; \omega$ with
\begin{align}
\label{E1.7}
\omega_+^\prime &= \lb 0; m+1, a_1, a_2, \ldots \rb  \\
\nonumber
\omega^\prime_- &= \lbd 0;m,b_1,b_2,\ldots \rbd  \in \left[r_q,-r_q\right].
\end{align}
Then the  bi-infinite sequence $\underline{a}^\prime:= \lbs \ldots,-b_2,-b_1,-m;m+1,a_1, a_2,\ldots \rbs$ is $q$-regular and by Lemma~\ref{D1.7} $\big(S\, \omega_+^\prime, - \omega_-^\prime\big) = \Pi(\underline{a}^\prime) \in \Omega_q$.

\smallskip

The case $\omega_+<0$ and $\omega_-\in \left[r_q,R_q\right] \subset I_{R_q}$ can be treated in the same way.

\medskip

The proof for $q = 3$ is similar to the case $q \geq 4$, however there are
the four cases $a_0 \geq 3$, $a_0 = 2$, $a_0 = 1$ and $a_0 = 0$ to be considered.

If $a_0 \geq 3$ then we can argue as in the case $a_0 \geq 2$ before.
Since $\omega_- < -r_3 = \lbd 0; \overline{-3} \rbd$, the bi-infinite sequence $\underline{a}$ in \eqref{E1.4} is $q$-regular and $(S\omega_+,-\omega_-) = (x,y) = \Pi(\underline{a}) \in \Omega_q$.

The cases $a_0 = 2$ and $a_0 = 1$ are similar to the case $a_0=1$ for $q \geq 4$: we just take the integer $\abs{m} \geq 5$ with $m \, b_1 <0$.
Then $\omega^\prime$ is defined as $\big(\omega^\prime_-,\omega^\prime_+\big) = \omega^\prime := g \; \omega$ with $g := ST^m$ and hence
\begin{align*}
\omega_+^\prime &=
\begin{cases}
\lb 0; m+2, a_1, a_2, \ldots \rb & \text{if $a_0=2$ and} \\
\lb 0; m+1, a_1, a_2, \ldots \rb & \text{if $a_0=1$},
\end{cases}  \\
\omega^\prime_- &= \lbd 0;m,b_1,b_2,\ldots \rbd  \in \left[r_3,-r_3\right].
\end{align*}

The case $a_0 = 0$ is similar to the case $a_0=0$ for $q \geq 4$, if we choose there the integer $m  \geq 3$ and recall $r_3 = \lb 0; \overline{3} \rb$ in \eqref{B3.8}.
\end{proof}

\subsection{The transfer operator for $\Gn{q}$}
\label{tr}
The authors of \cite{MS08} have constructed a Poincar\'e section $\Sigma$ for the geodesic flow $\Phi_t \colon
S_1 \,\Gn{q}\backslash\HH \to S_1 \,\Gn{q}\backslash\HH$ on the Hecke surfaces $\Gn{q} \backslash \HH$ for which the Poincar\'e map $P \colon \Sigma\to\Sigma$ is basically given by the natural extension $F_q$ of the map $f_q \colon \Iq \to \Iq$.
The periodic orbits of this geodesic flow can therefore be characterized by the periodic orbits of $F_q$ and therefore also by the periodic orbits of the map $f_q$ respectively its periodic points which determine the ones of $F_q$ uniquely.
Indeed, Theorem~\ref{B5.2} implies an almost one-to-one correspondence between the periodic orbits of the geodesic flow on the Hecke surfaces $\Gn{q}\backslash \HH$ and the periodic orbits of the map $f_q$, only the periodic orbits of the points $r_q$ and $-r_q$ which are not equivalent under the map $f_q$ lead to the same periodic orbit of the geodesic flow since these points are $\Gn{q}$-equivalent.
This shows already that the Selberg zeta function $Z_{\Gn{q}}$ for the Hecke triangle groups defined as
\[
Z_{\Gn{q}}(\beta) = \prod_{k=0}^\infty \prod _{\gamma \text{ prime}} \left( 1-e^{ -(\beta+k) l(\gamma)} \right),
\]
where the product is over the prime periodic orbits $\gamma$ of the geodesic flow and $l(\gamma)$ denotes its period (and hence the length of the corresponding closed geodesic), cannot be expressed in terms of the transfer operator for the map $f_q$ alone.
Indeed, to relate the above Selberg zeta function to the Poincar\'e map $P$ one uses the following Lemma by Ruelle \cite{R94}:

\begin{lemma}
\label{tr0}
$\displaystyle Z_{\Gn{q}}(\beta)= \prod_{k=0}^\infty e^{-\sum_{n=1}^\infty \frac{1}{n} Z_n(\beta +k)}$ for $\re{\beta} > 1$ where $Z_n(\beta)$ is the so called partition function $Z_n(\beta)=\sum_{x\in \mathrm{Fix}P^n}
e^{-\beta\sum_{l=0}^\infty r\big(P^l(x)\big)}$ and $r \colon \Sigma\to \R^+$ denotes the recurrence time of the geodesic flow with respect to the Poincar\'e section $\Sigma$.
\end{lemma}

In the transfer operator approach to the dynamical zeta functions the partition functions $Z_n(\beta)$ get expressed in terms of the traces of an operator constructed from the Poincar\'e map $P \colon \Sigma\to \Sigma$ respectively its restriction to the unstable directions.
In our case the unstable direction is one-dimensional and the restriction of $P$ to it is basically just the map $f_q:\Iq\to \Iq$.
On the other hand one knows that the recurrence time $r \colon \Sigma\to \R^+$ in our case is given by $r(x)=\log \abs{f_q^\prime(x)}$.
The Ruelle transfer operator $\cL_\beta$ then has  the following form
\begin{equation}
\label{tr1}
\cL_\beta g(x) = \sum_{y\in f_q^{-1}(x)} e^{-\beta \,r(y)} \, g(y)
\end{equation}
where $g \colon \Iq\to \mathbb{C}$ is some complex valued function and $\re{\beta} >1$ to ensure convergence of the series.
To get an explicit form for the operator $\cL_\beta$ one has to determine the preimages $y$ of any point $x \in \Iq$.
For this recall the Markov partition $\Iq = \bigcup_{i\in A_{\kappa_q}}\Phi_i$ with $A_{\kappa_q}=\{\pm 1,\ldots,\pm \kappa_q\}$ in (\ref{C3.6}), determined by the intervals $\Phi_i$, and the local inverses $\vartheta_{\pm m}(x) := \left(f_q\big|_{J_{\pm m}}\right)^{-1}(x) = \frac{-1}{x \pm m \lambda_q}$ on the intervals $J_{\pm m}$, $1\leq m\leq \infty$, respectively $2\leq m\leq \infty$ for $q=3$, defined in \S\ref{C3}.
For $1\leq i\leq \kappa_q$ denote by $\cN_i$ the set $\cN_i := \left\{n\in\Z\smallsetminus\{0\} \text{ such that there exists } j\in A_{\kappa_q} \text{ with } \vartheta_n (\Phi_i)\subset \Phi_j \right\}$.
But then  $\cN_i= \bigcup_{j\in A_{\kappa_q}}\cN_{i,j}$ with $\cN_{i,j}:= \left\{n \in \Z \smallsetminus\{0\} \text{ such that } \vartheta_n( \Phi_i)\subset \Phi_j\right\}$.
Using these sets we can rewrite the transfer operator $\cL_\beta$ in \eqref{tr1} as
\begin{equation}
\label{tr2}
\cL_\beta g(x)= \sum_{i\in A_{\kappa_q}} \chi_{\Phi_i}(x) \sum_{n\in\cN_i} \left(\vartheta_n^\prime(x)\right)^\beta \, g(\vartheta_n(x)),
\end{equation}
with $\chi_{\Phi_i}$ the characteristic function of the set $\Phi_i$.
With $g_i:= g\big|_{\Phi_i}$ this can be written also as follows
\begin{equation}
\label{tr3}
(\cL_\beta g)_i(x) = \sum_{j\in A_{\kappa_q} }\sum_{n\in\cN_{i,j}} \left(\vartheta_n^\prime(x)\right)^\beta \, g_j(\vartheta_n(x)), \quad x \in I_i.
\end{equation}
Thereby we used the Markov property of the partition $\Iq=\bigcup_{i\in A_{\kappa_q}}\Phi_i$.
If $g_i$ is continuous on $\Phi_i$ for all $i\in A_{\kappa_q}$ then also $(\cL_\beta g)_i$ is continuous on $\Phi_i$, that means $\cL_\beta$ maps piecewise continuous functions to piecewise continuous functions.
Unfortunately on the Banach space $B=\oplus_{i\in A_{\kappa_q}} C(\Phi_i)$ of piecewise continuous functions the operator $\cL_\beta$ is not trace class, it is even not compact.
Much better spectral properties however can be achieved by defining $\cL_\beta$ on a space of piecewise holomorphic functions.
This is possible since all the maps $\vartheta_{\pm m}, \, m\geq 1$ have holomorphic extensions to a complex neighbourhood of $\Iq$.
Indeed one shows
\begin{lemma}
\label{tr4}
There exist open discs $D_i\subset \C$, $i\in A_{\kappa_q} $ with $\Phi_i\subset D_i$ such that for all $n\in \cN_{i,j}$ we have $\vartheta_n(\overline{D_i})\subset D_j$.
\end{lemma}

Consider therefore the Banach space $B=\oplus_{i\in A_{\kappa_q}} B(D_i)$ with $B(D_i)$ the Banach space of holomorphic functions on the disc $D_i$ with the sup norm.
On this space the transfer operator $\cL_\beta$ has the form
\begin{equation}
\label{tr5}
(\cL_\beta g)_i(z)
=
\sum_{j\in A_{\kappa_q} } \sum_{n\in\cN_{i,j}}
   \left(\frac{1}{z+n\lambda_q}\right)^{2\beta} \, g_j \left(\frac{-1}{z+n\lambda_q} \right),
\qquad z\in D_i
\end{equation}
which is well defined for $\re{\beta} >\frac{1}{2}$.
In a forthcoming paper we will discuss the spectral properties of this operator and its relation to the Selberg zeta function for the Hecke triangle groups $\Gn{q}$.
Here we give the explicit form of this operator for the case $q=3$ and $q=4$.

For $q=3$ one has $\kappa_3=1$ and therefore $A_{\kappa_3}=\{\pm 1\}$.
The index sets $\cN_{i,j},\, i,j\in A_{\kappa_3} $ are given by $\cN_{1,1}=\Z_{\geq 3}$, $\cN_{1,-1} = \Z_{\leq -2}$, $\cN_{-1,1} = \Z_{\geq 2}$ and $\cN_{-1,-1}=\Z_{\leq -3}$.

For $q=4$ one has also $\kappa_4=1$ and hence $A_{\kappa_4}=\{\pm 1\}$.
The index sets $\cN_{i,j},\, i,j\in A_{\kappa_4} $ are given by $\cN_{1,1}=\Z_{\geq 2}$, $\cN_{1,-1}=\Z_{\leq -1}$, $\cN_{-1,1}=\Z_{\geq 1}$ and $\cN_{-1,-1}=\Z_{\leq -2}$.
This leads in these two cases to the following transfer operators
\begin{equation}
\label{tr6}
\begin{split}
(\cL_\beta g)_1(z)
=&
\sum_{n=3 \,(2)}^\infty \left(\frac{1}{z+n\lambda_q}\right)^{2\beta} g_1   \left(\frac{-1}{z+n\lambda_q }\right)\\
&+\sum_{n=2 \,(1)}^\infty \left(\frac{1}{z-n\lambda_q}\right)^{2\beta} g_{-1}\left(\frac{-1}{z-n\lambda_q}\right),
\quad z \in D_1,\\
(\cL_\beta g)_{-1}(z)
=&
\sum_{n=2 \,(1)}^\infty \left(\frac{1}{z+n\lambda_q}\right)^{2\beta} g_1   \left(\frac{-1}{z+n\lambda_q }\right)\\
&+\sum_{n=3 \,(2)}^\infty \left(\frac{1}{z-n\lambda_q}\right)^{2\beta} g_{-1}\left(\frac{-1}{z-n\lambda_q }\right),
\quad z\in D_2
\end{split}
\end{equation}
where $\lambda_3 =1$ and $\lambda_4=\sqrt{2}$ and the summation index in brackets belongs to the case $q=4$.
For $q=3,4$ the discs $D_i,\, i=\pm 1$ can be taken as $D_{\pm 1} = \pm \left\{z\in\C; \; \abs{z -\left(\frac{\lambda_q-2}{4}\right)} < \frac{\lambda_q+2}{4}\right\}$.

For $q=3$ this operator and its eigenfunctions with eigenvalue $\rho =1$ have been discussed in \cite{BM09} where it was shown that these eigenfunctions are directly related to the eigenfunctions with eigenvalues $\rho=\pm 1$ of the transfer operator for the modular group $\Gn{3}$ derived from a symbolic dynamics for the geodesic flow using the Gauss continued fractions in \cite{Ma91}.

\section{$\lambda_q$-CF's and Rosen $\lambda$-fractions}
\label{F}

\subsection{Regular $\lambda_q$-CF's and reduced Rosen $\lambda$-fractions ($q \geq 4$)}
\label{F1}
In \cite{Ro54} Rosen discussed continued fractions of the form
\begin{equation}
\label{F1.0}
[r_0;(\eps_1:r_1),(\eps_2:r_2),(\eps_3:r_3),\ldots]
=
r_0\lambda_q + \frac{\eps_1}{r_1\lambda_q + \frac{\eps_2}{r_2\lambda_q + \frac{\eps_3}{r_3\lambda_q + \ldots}}}
\end{equation}
with $r_0 \in \Z$ and $\eps_i = \pm 1$, $r_i \geq 1$ for $i \in \N$.
We call such expansions \emph{Rosen $\lambda$-fraction}.

Rosen $\lambda$-fractions and $\lambda_q$-CF's can easily be transformed into each other using the relations
\begin{align}
\label{F1.1}
& [r_0;(\eps_1:r_1),(\eps_2:r_2),(\eps_3:r_3),\ldots]  \\
\nonumber
&
= r_0\lambda_q + \frac{\eps_1}{r_1\lambda_q + \frac{\eps_2}{r_2\lambda_q + \frac{\eps_3}{r_3\lambda_q + \ldots}}}
= r_0\lambda_q + \frac{-1}{-\eps_1 r_1\lambda_q + \frac{-1}{\eps_1 \eps_2 r_2\lambda_q + \frac{-1}{-\eps_1 \eps_2 \eps_3 r_3\lambda_q + \ldots } } }  \\
\nonumber
&
= [ r_0;-\eps_1 \, r_1, \eps_1\eps_2 \, r_2, -\eps_1\eps_2\eps_3 \, r_3, \,\ldots\,  , (-1)^{i} \, \eps_1 \cdots \eps_i \, r_i, \, \ldots ]
\end{align}
and
\begin{align}
\label{F1.2}
&[ a_0;a_1,a_2,a_3,\ldots ]
= a_0\lambda_q + \frac{-1}{a_1\lambda_q + \frac{-1}{a_2\lambda_q + \frac{-1}{a_3\lambda_q + \ldots}}} \\
\nonumber
&\quad = a_0\lambda_q + \frac{-\sign{a_1}}{\abs{a_1}\lambda_q + \frac{-\sign{a_1}\sign{a_2}}{a_2\lambda_q + \frac{-\sign{a_2}\sign{a_3}}{a_3\lambda_q + \ldots } } }  \\
\nonumber
& \quad = [ a_0; (-\sign{a_1} : \abs{a_1}), (-\sign{a_1}\sign{a_2}: \abs{a_2}),  \\
\nonumber
& \quad \qquad\qquad\qquad\qquad\qquad\qquad     (-\sign{a_2}\sign{a_3}:\abs{a_3}) , \ldots ].
\end{align}

As claimed in \cite[Remark~15]{MS08} these relations imply directly
\begin{lemma}
\label{F1.3}
Given a Rosen $\lambda$-fraction $[r_0;\ldots,(\eps_i:r_i),(\eps_{i+1}:r_{i+1}),\ldots]$ and its corresponding $\lambda_q$-CF $[r_0 ;\ldots, a_i,a_{i+1}, \ldots]$ in \eqref{F1.2} we have
\[
\eps_1 = -\sign{a_1}
\quad \text{and} \quad
\eps_{i+1} = -\sign{a_i \, a_{i+1}}, \quad (i \in \N).
\]
\end{lemma}

\begin{proof}
Equation~\ref{F1.1} shows that the sign of the $i^\text{th}$ digit $(-1)^{i} \, \eps_1 \cdots \eps_i \, r_i$ in the formal $\lambda_q$-CF is determined by $(-1)^i\, \eps_1 \cdots \eps_i$.
Hence, the ratio
\[
\frac{(-1)^{i+1}\, \eps_1 \cdots \eps_i\eps_{i+1}}{(-1)^i\, \eps_1 \cdots \eps_i} = -\eps_{i+1}
\]
determines whether $a_i$ and $a_{i+1}$ have the same or opposite signs.
\end{proof}

\smallskip

Equations \eqref{F1.1} and \eqref{F1.2} indeed relate  regular $\lambda_q$-CF's and reduced Rosen $\lambda$-fractions as we show next.
Set
\begin{equation}
\label{F1.4}
h_\rR=h_\text{Rosen} := \nextinteger{\frac{q-3}{2}} =
\begin{cases}
h_q-1 & \text{if } q \text{ is even and} \\
h_q & \text{if } q \text{ is odd.}
\end{cases}
\end{equation}
The   above $h_\rR$ coincides with ``$h$'' used in \cite[above Definition~1]{Ro54}.
Recall from \cite[Page~555]{Ro54}
\begin{definition}[Reduced Rosen $\lambda$-fractions]
\label{F1.5}
The Rosen $\lambda$-fraction in \eqref{F1.0} is called \emph{reduced} if it satisfies the following conditions:
\begin{enumerate}
\item
\label{F1.5a}
Blocks of the form
\[
(\ast:1),\underbrace{(-1:1) , \ldots, (-1,1)}_{h_R \text{ times}} ,(-1:\ast)
\]
do not appear.
\item
\label{F1.5b}
For $q$ odd, blocks of the form
\[
( \ast:1),\underbrace{(-1:1),\ldots,(-1:1)}_{h_\rR \text{ times}}
\]
do not appear.
\item
\label{F1.5c}
For $q$ odd, blocks of the form
\[
( \ast:1),\underbrace{(-1:1),\ldots,(-1:1)}_{h_\rR - 1 \text{ times}}, ( -1:2),\underbrace{(-1:1),\ldots,(-1:1)}_{h_\rR \text{ times}},(-1,\ast)
\]
do not appear.
\item
\label{F1.5d}
For $q$ odd, a finite Rosen $\lambda$-fraction expansion does not terminate in a block of the form
\[
( \ast:1),\underbrace{(-1:1),\ldots,(-1:1)}_{h_\rR \text{ times}}.
\]
\item
\label{F1.5e}
The value $\pm \frac{\lambda_q}{2}$ of the tail $[(\ast:r_i),(\eps_{i+1}:r_{i+1}), \ldots, (\eps_{i+k}:r_{i+k})]$ of a finite Rosen $\lambda$-fractions leads because of $r_{i-1}\lambda_q \pm \frac{\lambda_q}{2} = (r_{i-1} \pm 1)\lambda_q \mp \frac{\lambda_q}{2}$  to non-uniqueness of the expansion.
We always choose the first possibility.
\end{enumerate}
\end{definition}

Then one shows
\begin{lemma}
\label{F1.6}
The following three statements hold:
\begin{itemize}
\item
The $\lambda_q$-CF associated to a reduced Rosen $\lambda$-fraction expansion in \eqref{F1.1} is regular.
\item
The Rosen $\lambda$-fraction expansion corresponding to a regular $\lambda_q$-CF in \eqref{F1.2} satisfies Properties (\ref{F1.5a})--(\ref{F1.5d}) of Definition~\ref{F1.5}.
\item
The two expansions of the finite Rosen $\lambda$-fractions in (\ref{F1.5e}) of Definition~\ref{F1.5} correspond to the identities of the finite regular $\lambda_q$-CF's in Lemma~\ref{B2.7}.
\end{itemize}
\end{lemma}

\begin{proof}
Let $x\in \R$ have the regular Rosen $\lambda$-fraction expansion~\eqref{F1.0}.
We have to show that the corresponding $\lambda_q$-CF in~\eqref{F1.1} does not contain any forbidden block from $\cB_q$.
We consider the cases $q$ even and $q$ odd separately.

Let $q$ be even.
Using Lemma~\ref{F1.3} and the identity $h_\rR=h_q-1$ in \eqref{F1.4} we see that Property~(\ref{F1.5a}) of Definition~\ref{F1.5} corresponds to the absence of blocks of the form $\big[(\pm 1)^{h_q},\pm m \big]$ for any $m \in \Z_{\geq 1}$.

Consider next $q$ odd.
Using again Lemma~\ref{F1.3} and the identity $h_\rR=h_q$ we see that Property~(\ref{F1.5b}) of Definition~\ref{F1.5} corresponds to the absence of blocks of the form $\big[(\pm 1)^{h_q+1} \big]$.
Similarly, Property~(\ref{F1.5c}) corresponds to the absence of blocks of the form $\big[(\pm 1)^{h_q},\pm 2, (\pm 1)^{h_q},\pm m \big]$ for any $m \in \Z_{\geq 1}$.

This shows that no forbidden block from $\cB_q$ appears in the $\lambda_q$-CF in \eqref{F1.1}.

\smallskip

Next, let $x \in \R$ have the regular $\lambda_q$-CF $x= \lb a_0;a_1,a_2,a_3,\ldots \rb$.
We have to show that the corresponding formal Rosen $\lambda$-fraction in~\eqref{F1.2} satisfies properties (\ref{F1.5a})--(\ref{F1.5d}) of Definition~\ref{F1.5}.
Again, we discuss the cases $q$ even and $q$ odd separately.

Consider first $q$ even.
Using Lemma~\ref{F1.3} and the identity $h_\rR=h_q-1$ in \eqref{F1.4} we find that forbidden blocks of the form $\big[(\pm 1)^{h_q},\pm m\big]$  for any $m \in \Z_{\geq 1}$ imply Property~(\ref{F1.5a}).
Property~\ref{F1.5e} corresponds just to the ambiguity of finite $\lambda_q$-CF's  given in Lemma~\ref{B2.7} since the tails $\big[(\pm 1)^{h_q}\big]$  correspond to $\mp\frac{\lambda_q}{2}$.

Consider next $q$ odd.
Using Lemma~\ref{F1.3} and recalling the identity $h_\rR=h_q$ in \eqref{F1.4} we see that forbidden blocks $\big[(\pm 1)^{h_q+1}\big]$  imply Property~(\ref{F1.5b}) and Property~(\ref{F1.5d}) of Definition~\ref{F1.5}.
Property~(\ref{F1.5b}) also implies Property~(\ref{F1.5a}).
Property~(\ref{F1.5c}) follows from the forbidden blocks $\big[(\pm 1)^{h_q},\pm 2,(\pm 1)^{h_q},\pm m\big]$ and for $m \in \Z_{\geq 1}$.
The ambiguity of the regular $\lambda_q$-CF's in Lemma~\ref{B2.7} implies Property~(\ref{F1.5e}).

\smallskip

To finish the proof of Lemma~\ref{F1.6} consider $q$ even and the finite regular $\lambda_q$-CF $\lb a_0; a_1, \ldots, a_n,(1)^{h_q} \rb$.
Using Equation~\eqref{F1.2} we rewrite it as the Rosen $\lambda$-fraction
\begin{multline}
\nonumber
\big[ a_0; (-\sign{a_1} : \abs{a_1}), \ldots  \\
 \ldots, (-\sign{a_{n-1}}\sign{a_n}:\abs{a_n}) ,  (-\sign{a_n}:1), (-1:1)^{h_q-1}\big].
\end{multline}
Since by the equation following~(4.2) in \cite{Ro54} and by \cite[(4)]{BKS99} the identity $\big[0; \linebreak[3](1:1), (-1:1)^{h-1}\big] = \frac{\lambda_q}{2}$ holds, we are in the situation of Property~(\ref{F1.5e}).
If $a_n<0$ we choose the ``$+$''-sign in Property~(\ref{F1.5e}).
If $a_n>0$ we use Lemma~\ref{B2.7} to rewrite the finite regular $\lambda_q$-CF such that its tail ends in $\big[a_n-1,(-1)^{h_q}\big]$ with $\sign{a_n-1} = \sign{a_n}$.
Using Equation~\eqref{F1.2} we arrive at the Rosen $\lambda$-fraction
\begin{align*}
& \big[ a_0; (-\sign{a_1} : \abs{a_1}), \ldots  \\
& \quad \ldots, (-\sign{a_{n-1}}\sign{a_n-1}:\abs{a_n-1}) ,  (\sign{a_n-1}:1), (-1:1)^{h_q-1} \big].
\end{align*}
with the correct tail.

The case $q$ odd is analogous to $q$ even, with the only difference that the reduced Rosen $\lambda$-fraction $\frac{\lambda_q}{2} = \big[0;(1:1), (-1:1)^{h_q-1},(-1:2),(-1,1)^{h_q} \big]$ as given in \cite[(4)]{BKS99} has the corresponding tail $\big[(1)^{h_q},2,(1)^{h_q}\big]$.
\end{proof}

\begin{remark}
\label{F1.7}
Consider the $\lambda_q$-CF of $\pm r_q$ in \eqref{B3.8}.
Their corresponding Rosen $\lambda$-fractions according to formula~\eqref{F1.1} are
\[
r_q=
\begin{cases}
\big[ 0; (-1:1)^{h_q-1},\ov{(-1:2),(-1,1)^{h_q-1}} \big]    \qquad\qquad 	\text{for $q$ even  and} \\
\big[ 0; (-1,1),\ov{(-1:1)^{h_q-1},(-1:2),(-1:1)^{h_q-1},(-1:2),(-1,1)} \big]  \\
\qquad\qquad\qquad\qquad\qquad\qquad\qquad\qquad\qquad \;		 \text{for $q$ odd}
\end{cases}
\]
and
\[
-r_q=
\begin{cases}
[ 0; (1:2),\ov{(-1:2)} ]  \qquad\qquad\qquad\qquad\qquad\quad\;	\text{for } q =4, \\
[ 0; (1:1),(-1:1)^{h_q-2},\ov{(-1:2),(-1,1)^{h_q-1}}] \qquad 	\text{for even } q \geq 6,\\
[ 0; (1,1),\ov{(-1:1)^{h_q-1},(-1:2),(-1:1)^{h_q-1},(-1:2),(-1,1)} ]  \\
\qquad\qquad\qquad\qquad\qquad\qquad\qquad\qquad\qquad\quad 	\text{for $q$ odd},
\end{cases}
\]
where $(-1:1)^0$ means that the digit $(-1:1)$ is absent.
The Rosen $\lambda$-fractions hence have the same tail.
\end{remark}

\begin{remark}
\label{F1.8}
The generating map $f^\star_q$ for the dual regular $\lambda_q$-CF  and the generating map $f_q^\rR$ for the Rosen $\lambda$-fractions in  \cite{BKS99} satisfy
\[
f_q^\star(-x)
=
\frac{1}{x} - \lambda_q \nextinteger{\frac{1}{x \lambda_q + \frac{R_q}{2}}}
=
\frac{1}{x} - \lambda_q \nextinteger{\frac{1}{x \lambda_q + 1 + \frac{r_q}{2}}}
=
f_q^\rR(x)
\]
for all $x \in \left(0,\frac{\lambda_q}{2}\right)$.

Formally, we find also
\[
f_q^\star(-x)
=
\frac{1}{x} - \lambda_q \nextinteger{\frac{1}{x \lambda_q + \frac{R_q}{2}}}
=
\frac{1}{x} - \lambda_q \nextinteger{\frac{1}{x \lambda_q + 1 + \frac{r_q}{2}}}
=
T_{-\frac{r}{\lambda_q}}(x)
\]
for all $x \in \left(0,\frac{\lambda_q}{2}\right)$ where $T_\alpha$ is the generating map of the $\alpha$-Rosen fractions discussed in \cite{DKS09}.
However, the parameter $\alpha = -\frac{r_q}{\lambda_q}$ lies outside the range $\alpha \in \left[\frac{1}{2}, \frac{1}{\lambda_q} \right]$ discussed in \cite{DKS09}, since by \eqref{B3.8} and \eqref{B3.6}  $-r_q=\lambda_q-R_q \in \left( 0,\frac{\lambda_q}{2} \right)$.
\end{remark}

\subsection{Regular $\lambda_q$-CF's and convergents}
\label{F2}
We define the $n^\text{th}$ convergent in the following way:
\begin{definition}
\label{F2.1}
Given a regular $\lambda_q$-CF $\lb a_0; a_1, \ldots, a_n,\ldots \rb$ of length at least $n$ we define its \emph{$n^\text{th}$ convergent} as the fraction $\frac{p_n}{q_n}$ where the numerator $p_n$ and denominator $q_n$ are given as entries in the vector
\begin{equation}
\label{F2.2}
{p_n \choose q_n} = T_q^{a_0} \, S T_q^{a_1} \, S T_q^{a_2} \, \cdots \, S T_q^{a_n} \; {0 \choose 1}.
\end{equation}
\end{definition}

The convergents then satisfy the recursion relation
\begin{equation}
\label{F2.3}
{ p_n \choose q_n}
=
\Matrix{p_{n-2}}{p_{n-1}}{q_{n-2}}{q_{n-1}} \, ST_q^{a_n} \, { 0 \choose 1}
=
{a_n \lambda_q \, p_{n-1} - p_{n-2} \choose a_n \lambda_q \, q_{n-1} - q_{n-2}}
\end{equation}
which holds also for $n=0$ and $n=1$ if we define $p_{-1} = 1$, $p_{-2}=0$, $q_{-1}=0$ and $q_{-2}=-1$.

\begin{remark}
\label{F2.4}
In the case $q=3$ Definition~\ref{F2.1} of the $n^\mathrm{th}$ convergent coincides with the usual definition as the ratio $\frac{p_n}{q_n} = \lb a_0;a_1,a_2, \ldots, a_n \rb$, since
\[
\lb a_0;a_1,a_2, \ldots, a_n \rb = \Matrix{\ast}{p_n}{\ast}{q_n} \;0
\]
where the last expression is to be understood as a M\"obius transformation.
\end{remark}

The following lemma will show that regular $\lambda_q$-CF's are indeed well defined and determine real numbers.
This obviously is true for finite regular $\lambda_q$-CF's.

\begin{lemma}
\label{F2.6}
Let $[a_0;a_1,a_2,\ldots]$ be an infinite regular $\lambda_q$-CF  and denote by $\frac{p_n}{q_n}$ it's $n^\text{th}$ convergent.
Then for $q\geq 4$ the fraction $\frac{\sign{q_n} \, p_n}{\abs{q_n}}$ is the $n^\text{th}$ convergent of the corresponding reduced Rosen $\lambda$-fraction as defined in \cite[Definition~3]{Ro54};
for $q=3$ the fraction $\frac{p_n}{q_n}$ is a ``N\"aherungsbruch'' in the sense of Hurwitz \cite[\S2]{Hu89}.
\end{lemma}

\begin{proof}
The case $q=3$ has been shown in \cite{Hu89}.

Hence assume $q \geq 4$.
Since the regular $\lambda_q$-CF is infinite, we don't have the ambiguities in Lemma~\ref{F1.6}.
which shows that the corresponding Rosen $\lambda$-fraction
\[
[ a_0; (-\sign{a_1} : \abs{a_1}), (-\sign{a_1}\sign{a_2}: \abs{a_2}), \ldots]
\]
is reduced.
The $n^\text{th}$ convergent $\frac{P_n}{Q_n}$ of the reduced Rosen $\lambda$-expansion is well defined and satisfies $Q_n \geq 1$ by \cite[Lemma~4]{Ro54}.
We have
\begin{align*}
\frac{P_n}{Q_n}
&= [ a_0; (-\sign{a_1} : \abs{a_1}), (-\sign{a_1}\sign{a_2}: \abs{a_2}), \ldots, \\
& \qquad\qquad\qquad\qquad\qquad \ldots, (-\sign{a_{n-1}}\sign{a_n}: \abs{a_n})] \\
&= [ a_0; a_1, a_2, \ldots, a_n] = \frac{p_n}{q_n}
\end{align*}
with $p_n$ and $q_n$ satisfying \eqref{F2.2}.
Hence we find indeed $P_n = \sign{q_n} \, p_n$ and $Q_n = \abs{q_n}$.
\end{proof}

This lemma shows that the results on convergents in \cite{Ro54} hold also for the regular $\lambda_q$-CF's.
We collect the relevant results in \cite{Ro54} and \cite{Hu89} in the following
\begin{lemma}
\label{F2.7}
The convergents $\frac{p_n}{q_n}$ of an infinite regular $\lambda_q$-CF satisfy:
\begin{itemize}
\item $q_n \not = 0 $ and $\abs{q_n} \geq \abs{q_{n-1}}$. For $q=3$ we have $\abs{q_n} > \abs{q_{n-1}}$.
\item $\lim_{n \to \infty} \abs{q_n} \to \infty$.
\item The sequence $\left( \frac{p_n}{q_n} \right)_{n\in\N}$ is a Cauchy sequence.
\end{itemize}
\end{lemma}

\begin{proof}
This follows from Lemma~4, Lemma~5, Theorem~4 and the proof of Theorem~5 in \cite{Ro54} for $q \geq 4$ and for $q=3$ from \S2 and \S3 in \cite{Hu89}.
\end{proof}

\smallskip

Now we can (re-)define infinite regular $\lambda_q$-CF's in the following way:
Let $\lb a_0; a_1,a_2, \ldots \rb$ be a regular $\lambda_q$-CF.
We assign the value $x$ to the regular $\lambda_q$-CF expansion and write $x = \lb a_0; a_1,a_2, \ldots \rb$ where  $x$ is the limit of the  sequence of convergents (see Definition~\ref{F1.4} and Lemma~\ref{F2.7}) of the corresponding Rosen $\lambda$-fraction.

Then the following estimate for the approximation of $x$ by the convergents holds:
\begin{lemma}
\label{F2.9}
Let $\lb a_0; a_1,a_2, \ldots\rb$ be an infinite regular $\lambda_q$-CF  and denote its $n^\text{th}$ convergents by $\frac{p_n}{q_n}$.
There exists a constant $\kappa_q>0$, independent of $x$, such that
\[
\abs{x - \frac{p_n}{q_n}} \leq \frac{1}{\kappa_q q_n^2}
\]
holds for all $n$.
\end{lemma}

\begin{proof}
The lemma follows for $q \geq 4$ from Theorem~4.6 in \cite{BKS99} and for $q=3$ from Satz on page 383 in \cite{Hu89}.
\end{proof}

\begin{remark}
\label{F2.10}
Obviously, Lemma~\ref{F2.9} implies that infinite regular $\lambda_q$-CF's converge.
This gives another proof of part of Proposition~\ref{B2.5}.
\end{remark}

\section*{Acknowledgement}
The authors thank Fredrik Stroemberg for many helpful discussions.

%
%

\bibliographystyle{amsalpha}

\begin{thebibliography}{99999}
\raggedright

\bibitem{BKS99}
R.M.\ Burton, C.\ Kraaikamp and T.A.\ Schmidt,
\newblock \textit{Natural extensions for the Rosen fractions},
\newblock Transactions of the American Mathematical Society \textbf{352} (1999), 1277-1298. \newline
\newblock \href{http://dx.doi.org/10.1090/S0002-9947-99-02442-3}{\texttt{doi:10.1090/S0002-9947-99-02442-3}}

\bibitem{BLZ}
R.W.\ Bruggeman, J.\ Lewis and D.\ Zagier,
\newblock \textit{Period functions for Maass wave forms. II: Cohomology},
\newblock in preparation.

\bibitem{BM09}
R.W.\ Bruggeman and T.\ M\"uhlenbruch,
\newblock \textit{Eigenfunctions of transfer operators and cohomology},
\newblock J.\ of Number Theory \textbf{129} (2009), 158--181. \newline
\newblock \href{http://dx.doi.org/10.1016/j.jnt.2008.08.003}{\texttt{doi:10.1016/j.jnt.2008.08.003}}

\bibitem{CM99}
C.H.\ Chang and D.\ Mayer,
\newblock \textit{The transfer operator approach to Selberg's zeta function and modular Maass wave forms for $\mathrm{PSL}(2,\Z)$},
\newblock In Emerging Applications of Number Theory, eds.\ D.\ Hejhal, M.\ Gutzwiller et al., IMA Volumes \textbf{109} (1999), 72--142, Springer-Verlag New York, 1999. 

\bibitem{DKS09}
K.~Dajani, C.~Kraaikamp and W.~Steiner
\newblock \textit{Metrical theory for $\alpha$-Rosen fractions},
\newblock To appear in J.\ of the European Math.\ Soc.\ (2009). \newline
\newblock \href{http://arxiv.org/abs/math/0702516v2}{\texttt{arXiv:math/0702516v2 [math.NT]}}

\bibitem{He83}
E.\ Hecke,
\newblock \textit{Mathematische Werke},
\newblock Vandenhoeck \& Ruprecht, G\"ottingen, 1959.

\bibitem{Hej83}
D.A.\ Hejhal,
\newblock \textit{The Selberg trace formula for ${\rm PSL}(2,\,R)$}, Vol.\ 2,
\newblock Lecture Notes in Mathematics \textbf{1001}, Springer-Verlag, Berlin, 1983. 

\bibitem{Hu89}
A.\ Hurwitz,
\newblock \textit{\"Uber eine besondere Art der Kettenbruch-Entwicklung reeller Gr\"ossen},
\newblock Acta Mathematica \textbf{12} (1889), 367--405. \newline
\newblock \href{http://dx.doi.org/10.1007/BF02391885}{\texttt{doi:10.1007/BF02391885}}

\bibitem{KU07}
S.\ Katok,and I.\ Ugarcovic,
\newblock \textit{Symbolic dynamics for the modular surface and beyond},
\newblock Bulletin of the American Mathematical Society \textbf{44} (2007), 87--132. \newline
\newblock \href{http://dx.doi.org/10.1090/S0273-0979-06-01115-3}{\texttt{doi:10.1090/S0273-0979-06-01115-3}}

\bibitem{LZ01}
J.\ Lewis and D.\ Zagier,
\newblock \textit{Period functions for Maass wave forms. I},
\newblock Annals of Mathematics \textbf{153} (2001), 191--258.\\
\newblock \href{http://dx.doi.org/10.2307/2661374}{\texttt{doi:10.2307/2661374}}

\bibitem{Marklof03}
J.\ Marklof,
\newblock \textit{Selberg's trace formula: an introduction},
\newblock Proceedings of the International School ``Quantum Chaos on Hyperbolic Manifolds'' (Schloss Reisensburg, G\"unzburg, Germany, 4-11 October 2003),
\newblock To appear in \emph{Lecture Notes in Physics}, Springer-Verlag. \newline
\newblock \href{http://arxiv.org/abs/math/0407288}{\texttt{http://arxiv.org/abs/math/0407288}}

\bibitem{Ma91}
D.\ Mayer,
\newblock \textit{The thermodynamic formalism approach to Selberg's zeta function for $\mathrm{PSL}(2,\Z)$},
\newblock  Bull.\ Amer.\ Math.\ Soc.\ (N.S.) \textbf{25} (1991), 55--60. \newline
\newblock \href{http://dx.doi.org/10.1090/S0273-0979-1991-16023-4}{\texttt{doi:10.1090/S0273-0979-1991-16023-4}}

\bibitem{Ma}
D.\ Mayer,
\newblock \textit{Transfer operators, the Selberg-zeta function and Lewis-Zagier theory of period functions},
\newblock Lecture notes of a course given in G\"unzburg, Germany, 4-11 October 2003.
\newblock To appear in \emph{Lecture Notes in Physics}, Springer-Verlag. \newline
\newblock \href{http://www.dynamik.tu-clausthal.de/research/preprints/Mayer.001.ps}{\texttt{\scriptsize http://www.dynamik.tu-clausthal.de/research/preprints/Mayer.001.ps}}

\bibitem{MS08}
D.\ Mayer and F.\ Str\"omberg,
\newblock \textit{Symbolic dynamics for the Geodesic flow on Hecke surfaces},
\newblock Journal of Modern Dynamics \textbf{2} (2008), 581--627. \newline
\newblock \href{http://dx.doi.org/10.3934/jmd.2008.2.581}{\texttt{doi:10.3934/jmd.2008.2.581}}

\bibitem{Na95}
H.\ Nakada,
\newblock \textit{Continued fractions, geodesic flows and Ford circles},
\newblock in \textit{Algorithms, Fractals, and Dynamics}, Edited by T.\ Takahashi, Plenum Press, New York, 1995. \newline
\newblock \href{http://www.ams.org/mathscinet-getitem?mr=1402490}{\texttt{MR1402490}}

\bibitem{Ru08}
Z.\ Rudnick,
\newblock \textit{What is$\dots$ quantum chaos?},
\newblock Notices Amer.\ Math. Soc.\ \textbf{55} (2008), no.\ 1, 32--34. \newline
\newblock \href{http://www.ams.org/notices/200801/tx080100032p.pdf}{\texttt{http://www.ams.org/notices/200801/tx080100032p.pdf}}

\bibitem{R94}
D.\ Ruelle,
\newblock \textit{Zeta Functions for Piecewise Monotone Maps of the Interval},
\newblock CRM Monograph Series, \textbf{vol.~4} ,
\newblock AMS, Providence, Rhode Island, 1994. 

\bibitem{R04}
D.\ Ruelle,
\newblock \textit{Thermodynamic Formalism},
\newblock  Cambridge Mathematical Library
\newblock  Cmabridge University Press, Cambridge CB2 2RU, UK, 2004. 

\bibitem{Sa95}
P.\ Sarnak,
\newblock \textit{Arithmetic quantum chaos},
\newblock The Schur lectures (1992) (Tel Aviv), 183--236, Israel Math.\ Conf. Proc. \textbf{8}, Bar-Ilan Univ., Ramat Gan, 1995. \newline
\newblock \href{http://www.math.princeton.edu/sarnak/Arithmetic Quantum Chaos.pdf}{\texttt{\scriptsize http://www.math.princeton.edu/sarnak/Arithmetic Quantum Chaos.pdf}}

\bibitem{SS95}
T.A.\ Schmidt and M.\ Sheingorn,
\newblock \textit{Length spectra of the Hecke triangle groups},
\newblock Mathematische Zeitschrift \textbf{220} (1995), 369--397. \newline
\newblock \href{http://dx.doi.org/10.1007/BF02572621}{\texttt{doi:10.1007/BF02572621}}


\bibitem{Ro54}
D.\ Rosen,
\newblock \textit{A class of continued fractions associated with certain properly discontinuous groups},
\newblock Duke Mathematical Journal \textbf{21} (1954), 549--563. \newline
\newblock \href{http://dx.doi.org/10.1215/S0012-7094-54-02154-7}{\texttt{doi:10.1215/S0012-7094-54-02154-7}}


\end{thebibliography}


\end{document}